\documentclass[12pt]{amsart}
\usepackage{amsfonts, amssymb, amsmath, amsthm, eucal, latexsym, array}
\usepackage[colorlinks=true,linkcolor=blue,urlcolor=violet,citecolor=magenta]{hyperref}
\usepackage{color}
\usepackage{fullpage}
\usepackage{verbatim}

\input{cyracc.def}

\font\tencyr=wncyr10 %scaled \magstephalf
\font\tencyi=wncyi10 %scaled \magstephalf
\def\rus{\tencyr\cyracc}
\def\rusi{\tencyi\cyracc}

\numberwithin{equation}{section}

\newtheorem{thm}{Theorem}[section]
\newtheorem{lm}[thm]{Lemma} %%[section]

\newtheorem{prop}[thm]{Proposition}
 
\theoremstyle{remark}
\newtheorem{ex}[thm]{Example}%[section]
\newtheorem{rmk}[thm]{Remark}%[section]
\newtheorem{cnj}[thm]{Conjecture}

\newtheorem{alg}[thm]{Algorithm}

\theoremstyle{definition}
\newcommand{\gt}{\mathfrak}

\newcommand{\ff}{\mathbb F}
\newcommand{\SL}{{\rm SL}}
\newcommand{\GL}{{\rm GL}}

\newcommand{\Spin}{{\rm Spin}}

\newcommand{\ind}{{\rm ind\,}}
\newcommand{\diag}{{\rm diag}}

\newcommand{\Lie}{\mathrm{Lie\,}}

\newcommand{\Hom}{\mathrm{Hom}}

\newcommand{\ad}{\mathrm{ad\,}}

\newcommand{\g}{\mathfrak g}

\newcommand{\esi}{\varepsilon}

\newcommand{\VV}{{\mathbb V}}
\newcommand{\W}{{\mathbb W}}
\newcommand{\BZ}{{\mathbb Z}}

\newcommand{\GR}[2]{{\textrm{{\bf #1}}}_{#2}}

\renewcommand{\le}{\leqslant}
\renewcommand{\ge}{\geqslant}

\newfam\Bbbfam\newfam\eufam\newfam\eusfam%
\font\Bbbfont=msbm10 scaled 1200%
\font\Bbbsmallfont=msbm8%
\textfont\Bbbfam=\Bbbfont\scriptfont\Bbbfam=\Bbbsmallfont%%
\font\euzw=eufm10 scaled 1200%
\font\euac=eufm7 scaled 1200%
\font\euacc=eufm7 scaled 1000%
\textfont\eufam=\euzw\scriptfont\eufam=\euac%
\scriptscriptfont\eufam=\euacc%
\font\euszw=eusm10 scaled 1200%
\font\eusac=eusm7 scaled 1200%
\font\eusacc=eusm7 scaled 1000%
\textfont\eusfam=\euszw\scriptfont\eusfam=\eusac%
\scriptscriptfont\eusfam=\eusacc%

\def\varnothing{\hbox {\Bbbfont\char'077}}

\newfam\eusfam%
\font\euszw=eusm10 scaled 1200%
\font\eusac=eusm7 scaled 1200%
\font\eusacc=eusm7 scaled 1000%
\textfont\eusfam=\euszw\scriptfont\eusfam=\eusac%
\scriptscriptfont\eusfam=\eusacc%
\usepackage{longtable}

\newcommand{\Z}{\mathbb{Z}}

\newcommand{\rank}{\mathrm{\mathop{rank}}}

\newcommand{\mf}[1]{\mathfrak{#1}}

\begin{document}
\hfill {\scriptsize March 22, 2010}%%% \today}
\vskip1ex

\title[GIB of $\theta$-representations]
{Good index behaviour of $\theta$-representations, I} 
\author{Willem A. de Graaf}
\address{Dipartimento di Matematica,
Universit\`{a} di Trento, via Sommariver 14, 
I-38100 Povo (Trento) Italy}
\email{degraaf@science.unitn.it}
\author{Oksana S. Yakimova}
\address{Mathematisches Institut, 
Universit\"at Erlangen-N\"urnberg,
Bismarckstrasse 1\,1/2,
91054 Erlangen  Germany}
\email{yakimova@mpim-bonn.mpg.de}
%\thanks{}
%\keywords{Nilpotent orbits, centralisers, symmetric invariants}
%\subjclass[2000]{17B45}
\begin{abstract}
Let $Q$ be an algebraic group with $\gt q=\Lie Q$ and $V$ a $Q$-module. The index of $V$ is the minimal codimension of the $Q$-orbits in the dual space $V^*$. There is a general inequality, due to Vinberg, relating the index of $V$ and the index of a $Q_v$-module $V/\gt q{\cdot}v$ for  $v\in V$. 
A pair $(Q,V)$ is said to have GIB if Vinberg's inequality turns into  
an equality for all $v\in V$. 
In this article, we are interested in the GIB property of $\theta$-representations,
where $\theta$ is a finite order automorphism of a simple Lie algebra $\gt g$. 
An automorphism of order $m$ defines a $\BZ/m\BZ$-grading 
$\gt g=\bigoplus\gt g_i$. If $G_0$ is the identity
component of $G^\theta$, then it acts on $\gt g_1$ 
and this action is called a $\theta$-representation.  
We classify inner automorphisms 
of $\gt{gl}_n$ and all finite order autmorphisms of the exceptional Lie algebras 
such that $(G_0,\gt g_1)$ has GIB and $\gt g_1$ contains a semisimple element. 
\end{abstract}
\maketitle

%%%\tableofcontents

%%%%%%%%%%%%%%%%%%%%%%%%%%%%%%%%
\section{Introduction}
%%%%%%%%%%%%%%%%%%%

Let $\mf{q}$ be a Lie algebra  over 
an algebraically closed field $\ff$ of characteristic zero 
and $V$ a finite-dimensional $\mf{q}$-module.
For  $\xi\in V^*$ we set
$$\mf{q}_\xi = \{ x\in \mf{q} \mid x{\cdot}\xi = 0\}.$$ 
Then the non-negative integer
$$\dim V - \max_{\xi \in V^*} (\dim \mf{q}{\cdot} \xi) = 
\dim V -\dim \mf{q} +\min_{\xi\in V^*}(\dim \mf{q}_\xi) $$
is called the {\it index} of $V$; and denoted by $\ind(\mf{q},V)$.

Suppose that $\mf{q}$ is the Lie algebra of an algebraic group $Q$,
and that $V$ is also a $Q$-module. 
Then by the Rosenlicht theorem, $\ind(\gt q,V)$ is equal to
${\rm tr.deg}\,\ff(V^*)^Q$.
Following \cite{py-gib} 
we say that the pair $(Q,V)$ %%% is defined
has a {\it good index behaviour\/} (GIB), if 
%%% to have GIB  
\begin{equation}\label{eqn:3} 
\ind( \mf{q}, V^*) = \ind(\mf{q}_v, (V/\mf{q}{\cdot} v)^*)
\end{equation}
for every $v\in V$. (Note that $V/\mf{q}{\cdot} v$ is a $\mf{q}_v$-module.) 
%%% It is known that 
It was noticed by Vinberg that
the left hand side is
less than or equal to the right hand side, see \cite[Sect.\,1]{Dima2}. 
%%% ; this inequality is due to Vinberg. 
Further, $(Q,V)$ is said to have GNIB (Good {\it Nilpotent} Index
Behaviour) if (\ref{eqn:3}) holds for all nilpotent elements $v\in V$ (where $v\in V$ 
is said to be nilpotent if $0\in \overline{Q{\cdot} v}$). 

If $v=0$, or $\dim (Q{\cdot}v)=\dim V-\ind(\gt q,V^*)$, 
or the stabiliser  $Q_v$ is reductive, then $v$ satisfies (\ref{eqn:3}),
see \cite{py-gib}. In general, it is a rather intricate problem to 
check the equality. One of the possible ways to prove that it holds 
for $v$
is to find $\xi\in V/\gt q{\cdot}v$ 
such that 
\begin{equation}\label{point}
\dim\gt q_v-\dim(\gt q_v)_\xi=\dim V-\dim\gt q{\cdot}v-\ind(\gt q,V).
\end{equation}
Note also that if (\ref{eqn:3}) is satisfied for $v$, then it is 
satisfied for all elements of the orbit $Q{\cdot}v$ as well.  

Checking GIB for a representation is even more complicated.
No general principle exists at the moment. The only method is 
to classify the $Q$-orbits and then compute the index for 
all of them.  A few positive results are known, for example, 
all representations of an algebraic torus $(\ff^{^\times})^m$ do have GIB 
(\cite{py-gib}). It would be interesting to understand what properties 
of a representation cause GIB to hold. 

Suppose that $V$ has only finitely many nilpotent 
$Q$-orbits. In that case the representation of $Q$ on 
$V$ is said to be {\it observable} and the GIB property is 
equivalent to GNIB, see \cite[Theorem~2.3]{py-gib}.  
Many observable representations %%% naturally 
arise in the context of reductive Lie algebras and their 
semisimple automorphisms. 

Let $G$ be a connected reductive algebraic group defined over $\ff$
and set $\gt g:=\Lie G$. The
group $G$ acts on $\g$ via the adjoint representation 
and this action is known to be observable.
In general, if $V=\gt q^*$, then 
%%% one writes simply 
$\ind(\gt q,\gt q)$ is the index,   
$\ind\gt q$, of $\gt q$ in the sense of Dixmier. If $\gamma\in\gt q^*$, 
then 
$\gt q^*/(\gt q{\cdot}\gamma)\cong\gt q_\gamma^*$ as a $Q_\gamma$-module.  
Therefore Vinberg's inequality reads $\ind\gt q\le\ind\gt q_\gamma$. 
Elashvili conjectured that all reductive Lie algebras 
have GIB. The conjecture is proved on a case-by-case basis 
in \cite{graaf} (exceptional Lie algebras) and
\cite{fan} (classical Lie algebras). 
An alternative proof is recently obtained by Charbonnel and Moreau
(\cite{ch-m}).

Let $\theta$ be an involution of $\gt g$. 
Then $\gt g=\gt g_0\oplus\gt g_1$, where $\g_i$ is the eigenspace 
of %%% an involution
$\theta$, %%% of $\g$, 
corresponding to the eigenvalue $(-1)^i$. Here $\gt g_0$ is a reductive 
subalgebra, which is the Lie algebra of the connected reductive
subgroup $G_0\subset G$. The subgroup $G_0$ acts on $\gt g_1$ via 
restriction of the adjoint representation of $G$.
In many aspects this representation is similar to the 
adjoint action of a reductive group. For example, 
all maximal subalgebras in $\gt g_1$ consisting of semisimple elements
are conjugate under $G_0$ (\cite{kr}). Such a subalgebra 
$\gt c\subset\gt g_1$ 
is usually referred to as a Cartan subspace and 
$\ind(\gt g_0,\gt g_1^*)=\dim\gt c=:\rank(G_0,\gt g_1)$. 
Kostant and Rallis (\cite{kr}) have shown also that the representation of 
$G_0$ on $\gt g_1$ is observable. 
As was found out in \cite{py-gib}, not all pairs $(G_0,\gt g_1)$ 
satisfy the GIB property. Therefore it is an 
interesting problem to describe those of  them, which do have GIB.  
In \cite{py-gib}, the GIB property was checked for all involutions 
except the following two (we give the corresponding symmetric pairs
$(\gt g,\gt g_0)$): $(\GR{E}{6},\gt{so}_{10}{\oplus}\ff)$, 
$(\GR{E}{7},\GR{E}{6}{\oplus}\ff)$. The calculations reported on in
this paper show that these two involutions do have GIB.

In the 70-s Vinberg (\cite{vinberg}) generalised results 
of Kostant and Rallis (\cite{kr}) to 
the set-up of arbitrary semisimple automorphisms of $\gt g$. 
Let $\theta$ be an automorphism of $\g$ of order $m$ and $\zeta$ 
is a primitive $m$-th root of unity. 
Then there is a $\Z/m\Z$-grading of $\g$,
$$ \g = \g_0 \oplus \g_1\oplus \cdots \oplus \g_{m-1},$$
where $\g_i$ is the eigenspace of $\g$ corresponding to $\zeta^i$.
Let $G_0\subset G$ be a connected algebraic group 
with the  Lie algebra $\g_0$. Then $G_0$ is reductive and it acts on $\g_1$ in a natural way. 
The group $G_0$, together with its action on $\g_1$, is called a $\theta$-group. 
The representation of $G_0$ on $\g_1$ is also called a $\theta$-representation.

Similar to the symmetric space situation a Cartan subspace of $\g_1$ is defined to be a maximal
subspace consisting of commuting semisimple elements. All Cartan
subspaces are conjugate under $G_0$, and 
the dimension of any of them 
%%% their dimension 
is called the rank of $\g_1$ 
(or rather of the pair $(G_0,\g_1)$, or of the 
$\theta$-representation afforded by $G_0$ and $\g_1$). 
According to \cite{vin-nilp}, \cite{vinberg}, all $\theta$-representations are 
observable. As a consequence, there is 
always a nilpotent orbit of dimension  
$\dim\gt g_1-\rank(G_0,\gt g_1)$.
In \cite{vin-nilp} Vinberg developed a method for classifying the nilpotent
$G_0$-orbits in $\g_1$.

In this paper we classify finite order automorphims 
of the exceptional Lie algebras and inner (finite order) automorphisms 
of $\gt{gl}_n$ such that $\rank(G_0,\gt g_1)>0$ and 
$(G_0,\gt g_1)$ has GIB. For the exceptional case 
the answer is given in Tables~\ref{tab:gibE6} to \ref{tab:gibG2} 
in Section~\ref{sec-ex}.
In the $\gt{gl}_n$ case a positive rank $\theta$-representation 
(with inner $\theta$) has GIB if and 
only if  
$\theta$ is a conjugation by one of the following diagonal matrices:
\begin{itemize}
\item  \ $\diag(1^{n-2},(-1)^2)$, $\diag(1^{n-1},-1)$, $\diag(1^3,(-1)^3)$ ($n=6$), 
\ see \cite{py-gib}; 
\item \ $\diag(1^2,\zeta^2,(\zeta^2)^{n-4})$, where $\zeta$ is a primitive 
$3$-d root of unity, \ 
see Theorem~\ref{m=3};
\item \ $\diag(1,\zeta^{r_1},\ldots,(\zeta^{m-1})^{r_{m-1}})$, where $\zeta$ 
is a primitive $m$-th root of unity and
there are no subsequences 
$r_i,r_{i+1},r_{i+2}$ with all elements 
being larger than $1$, 
see Theorem~\ref{4-2-1} and Proposition~\ref{1-groups-1}.
\end{itemize} 
In the $\gt{gl}_n$ case, the 
classification of nilpotent $G_0$-orbits, first obtained 
by Kempken~(\cite{Gisela-dis}), is presented in subsection~\ref{nilp-g1}.
In the exceptional case, we get representatives of the nilpotent 
orbits using the algorithms of \cite{gra15}. After that 
GIB is checked for each of them with Algorithm~\ref{alg_gib}.
The answer, automorphisms $\theta$ such that $(G_0,\gt g_1)$ 
has GIB, is given in terms of the so-called {\it Kac diagrams} 
(\cite{kac_autom}). Here
we summarise the main properties that we use of these diagrams; 
for a  more detailed explanation see \cite[Chapter X, \S5]{helgason}.

A Kac diagram is either an extended Dynkin diagram of $\gt g$ or the one 
obtained from it by gluing together points in the orbits of the diagram 
automorphism.  In addition, one attaches non-negative numbers, labels,
to the 
vertices. We are more interested in the $\theta$-representation, than 
in the $\theta$ itself and there is an easy way to read this from the 
Kac diagram. Assume for simplicity that $\theta$ is inner 
(in this paper $\theta$ is outer only for $\gt g=\GR{E}{6}$). 
Then $G_0$ contains a maximal torus of $G$ and the semisimple 
part of $\gt g_0$ is generated by all roots that have labels 
$0$ on the Kac diagram. The lowest weights of $\gt g_1$ (with 
respect to $G_0$) are in one-to-one correspondence 
with the roots labeled with $1$.  Finally, 
we notice that if $\rank(G_0,\gt g_1)>0$, then 
the Kac diagram has only labels  $0$ and $1$ (\cite{vinberg}). 
Hence we give these labels by colouring the nodes:
black means that the label is $1$, otherwise the label is $0$. 
For a diagram with all nodes black, $G_0$ is a torus and 
$(G_0,\g_1)$ has GIB by \cite[Proposition 1.3]{py-gib}. 
We do not put such ``all black'' diagrams in the tables. 

\vskip0.5ex

The GIB property of positive rank automorphisms in other 
classical types and outer automorphisms of $\gt{gl}_n$ 
will be studied in a forthcoming paper. 
Due to the large amount of cases and presumably a rather involved answer 
(cf. Proposition~\ref{prop}) we leave rank zero automorphisms 
aside. Inner automorphisms in type $A$ provide a remakable instance, where  
$\theta$-representations are also quiver representations. 
It would be interesting to check GIB for other (observable) 
quiver representations.

\vskip1ex

\noindent
{\bf Acknowledgements.}  
Parts of this work were carried out during our meetings in various 
mathematical institutions.  We are grateful to the 
University of Trento,  Centro Ennio De Giorgi (SNS, Pisa), and 
Emmy-Noether-Zentrum (Erlangen) for hospitality and  financial 
support. The second author would like to thank 
the Max-Planck-Institut f\"ur Mathematik (Bonn) for 
the excellent working conditions it provides.

\section{Some remarks on the GIB property}\label{alg-sec}

We begin this section
with a more explicit description of GIB for $\theta$-representations. 
Let $\gt g$ be a semisimple Lie algebra and 
$\gt g=\bigoplus\gt g_i$ a $\mathbb Z/m\mathbb Z$-grading.
For $x\in\gt g$, let $\gt g_x$ denote the centraliser of $x$ in $\gt g$
and set
%%% $\gt g_{0,x}$ denote the stabiliser of $x$ in $\gt g_0$ and
$\gt g_{i,x}:=\gt g_x\cap\gt g_i$. 
%%%%the stabiliser in $\gt g_i$.
For $x\in\gt g_1$, the inclusion $0\in\overline{G_0{\cdot}x}$ holds if and only if 
$x$ is a nilpotent element of $\gt g$ in the usual sense, 
see \cite{vinberg}.

\begin{prop}\label{prop:2}
%%%\begin{enumerate}
\begin{itemize}
\item[\sf(i)] \ 
%% \item 
For all nilpotent elements $e\in\g_1$ we have 
$$\rank(G_0,\g_1) \leq \ind(\g_{0,e},\g_{-1,e}).$$
\item[\sf(ii)] \ The pair $(G_0,\g_1)$ has GIB if and only if for all nilpotent 
elements $e\in\g_1$
we have $\rank(G_0,\g_1) = \ind(\g_{0,e},\g_{-1,e})$.
%%%\end{enumerate}
\end{itemize}
\end{prop}
\begin{proof} (cf. \cite[Proposition~2.6]{py-gib}). 
In this case (\ref{eqn:3}) translates to
$$\ind(\g_0,\g_1^*) = \ind(\g_{0,e}, (\g_1/[\g_0,e])^* ).$$
%%% $$\ind(\g_0,\g_1^*) \leq \ind(\g_{0,e}, (\g_1/[\g_0,e])^* ).$$
Now $\ind(\g_0,\g_1^*)$ is equal to $\dim \g_1 - \max_{x\in\g_1} \dim[\g_0,x]$.
By \cite[Theorem~5]{vinberg}, this is equal to $\rank( G_0, \g_1)$. 
The Killing form $\kappa$ gives a 
nondegenerate pairing $\kappa : \g_{-1}\times \g_1 \to \ff$. Using this 
pairing we get an isomorphism of $\g_0$-modules $\g_1^*\cong \g_{-1}$.
In the same way we get an isomorphism of $\g_{0,e}$-modules
$$(\g_1/[\g_0,e])^* \cong \{ y\in \g_{-1} \mid \kappa(y,[\g_0,e])=0 \}.$$
Let $y$ lie in the latter space. Then $0= \kappa( y,[\g_0,e]) = 
\kappa([e,y],\g_0)$. Now $[e,y]\in \g_0$ and $\kappa : \g_0\times \g_0
\to \ff$ is nondegenerate. Hence $[e,y]=0$, and it follows that
$\{ y\in \g_{-1} \mid \kappa(y,[\g_0,e])=0 \} = \g_{-1,e}$. 
Therefore $(\g_1/[\g_0,e])^*\cong\g_{-1,e}$ and  Vinberg's inequality
turns into   
$\rank(G_0,\g_1) \leq \ind(\g_{0,e},\g_{-1,e})$. 
Finally note that the $G_0$-module $\gt g_1$ has only finitely many
nilpotent orbits, see \cite{vinberg}, and thereby GIB is equivalent to GNIB
by \cite[Theorem~2.3]{py-gib}. 
\end{proof}

We will also need some easy technical statements concerning GIB.  

\begin{prop} Suppose that $\gt g=\Lie G$ is a reductive Lie algebra
and $V$ a finite dimensional $G$-module.
Then $\ind(\gt g,V)=\ind(\gt g,V^*)$, also $V$ and $V^*$ 
have (or do not have) GIB at the same time. 
\end{prop}
\begin{proof}
Each reductive group $G$ posses a so-called {\it Weyl involution}
$\sigma$, see e.g. \cite[Chapter~IX, \S5]{helgason}, which has a property that 
$\rho{\circ}\sigma\cong\rho^*$ for all finite dimensional representations
$\rho$. In particular, for each $v\in V$, where is a vector 
$\sigma(v)\in V^*$ such that $G_{\sigma(v)}=\sigma(G_v)$. 
Since $\sigma$ preserves the dimension of subgroups, 
we have   $\ind(\gt g,V)=\ind(\gt g,V^*)$. Moreover, $\sigma$ establishes 
an isomorphism between the representations of 
$G_v$ on $V/\gt g{\cdot}v$ and $G_{\sigma(v)}$ on $V^*/\gt g{\cdot}\sigma(v)$.
If GIB fails for $v$, it also fails for $\sigma(v)$.
\end{proof}

Since $\gt g_{{-}1}\!\cong\gt g_1^*$, these $\theta$-representations 
are of the same rank and GIB holds for 
$(G_0,\gt g_1)$ if and only if it holds for $(G_0,\gt g_{-1})$. 

\begin{prop}\label{summa}
A representation of $\gt q$ on $V$ has GIB if and only if 
for all $v\in V$ there is $w\in V$ such that 
$\dim(\gt q{\cdot}v+\gt q_v{\cdot}w)=\dim V-\ind(\gt q,V^*)$. 
\end{prop}
\begin{proof}
By definition, the representation has GIB if and only if 
for all $v\in V$ we have 
$\ind(\gt q_v,(V/\gt q{\cdot}v)^*)=\ind(\gt q,V^*)$.
The equality holds if and only if there is a coset 
$\bar w=w+\gt q{\cdot}v$ (with $w\in V$) 
such that 
$\dim(V/\gt q{\cdot}v)-\dim(\gt q_v{\cdot}\bar w)=\ind(\gt q,V^*)$. 
It remains to notice that the left hand side is equal to 
$\dim V-\dim(\gt q{\cdot}v)-\dim(\gt q_v{\cdot}\bar w)$ and 
that  
$\dim(\gt q{\cdot}v)+\dim(\gt q_v{\cdot}\bar w)=\dim(\gt q{\cdot}v+\gt q_v{\cdot}w)$.
\end{proof} 

Let $Q$ be an algebraic group acting on a finite dimensional 
vector space $V$ and $\gt q=\Lie Q$. 
Suppose that $\ind(\gt q,V^*)=1$ and
$av\not\in\gt q{\cdot}v$ 
for generic $v\in V$ and non-zero 
$a\in\ff$. Consider the action of 
$\widetilde{Q}:=Q{\times}\ff^{^\times}$ on $V$ such that
$t{\cdot}v=tv$ for all $t\in\ff^{^\times}$, $v\in V$. 
Set $\widetilde{\gt q}:=\Lie\widetilde{Q}$.
The two groups $Q$ and $\widetilde{Q}$ have 
different generic orbits on $V$. Hence 
$\ind(\widetilde{\gt q},V^*)=0$.  

\begin{lm}\label{ind-0-1} 
Suppose that the above assumptions are satisfied and 
the action of $\widetilde{Q}$ on $V$ has GIB. Then
the action of $Q$ on $V$ has GIB as well. 
\end{lm}
\begin{proof}
Take any $v\in V$. Then either 
$\gt q{\cdot} v=\widetilde{\gt q}{\cdot}v$ or 
$\dim \widetilde{\gt q}{\cdot}v=\dim \gt q{\cdot}v+1$ 
and  $\gt q_v=\widetilde{\gt q}_v$. If the first case takes place, 
then $V/\gt q{\cdot}v=V/\widetilde{\gt q}{\cdot}v$ and 
$\dim\gt q_v=\dim\widetilde{\gt q}_v-1$. Hence 
$\ind(\gt q_v, (V/\gt q{\cdot}v)^*)
 \le \ind(\widetilde{\gt q}_v,(V/\widetilde{\gt q}{\cdot}v)^*)+1=1$. 
In the second case 
$\gt q_v{\cdot}x+\widetilde{\gt q}{\cdot}v=V$ for  generic 
$x\in V$. Since $\dim \gt q{\cdot}v=\dim \widetilde{\gt q}{\cdot}v-1$, 
we again get the inequality  
$\ind(\gt q_v, (V/\gt q{\cdot}v)^*) \le 1$. 
\end{proof}

\begin{rmk}
The inverse implication is not true in general. 

Let
$\gt q=\left\{  \left(\begin{array}{ccc} 
                                      -t & s & 0 \\
                                       0 & t & 0 \\
                                       0 & 0 & t \\ 
                                       \end{array}\right) \mid t,s\in\mathbb F\right\}$
be a                                        
two-dimensional subalgebra  of $\gt{gl}_3$. Then
$\ind(\gt q,(\ff^3)^*)=1$ and the defining representation 
of $\gt q$ on $\ff^3$ has GIB. If we add a one dimensional central torus, 
then $\ind(\widetilde{\gt q},(\ff^3)^*)=0$, but the GIB property 
is not satisfied for the second basis vector. 
\end{rmk}

One of the ways to compute $\ind(\gt q,V)$ 
is related to the matrix $(\gt q{\cdot}V)$ of the action of 
$\gt q$ on $V$. Let $x_1,\ldots,x_n$ be
a basis of $\mf{q}$ and  $v_1,\ldots,v_s$ a basis of $V$. 
Then $(\gt q{\cdot}V)$ is an $n{\times}s$-matrix with 
entries $x_i{\cdot}v_j$. Each element of $V$ can be considered 
as a linear (or rational) function on $V^*$. 
Therefore it is possible to compute the rank of $(\gt q{\cdot}V)$
over a field $\ff(V^*)$.

\begin{lm}\label{matrix}
 We have $\ind(\gt q,V)=\dim V-\rank(\gt q{\cdot}V)$.
\end{lm}
\begin{proof}
Take $\xi\in V^*$ and
set $c_{ij}:=\xi(x_i{\cdot}v_j)$. 
Let 
$x=\sum_{i}\alpha_ix_i\in\gt q$ (with $\alpha_i\in\ff$).
Then
$$
x{\cdot}\xi(v_j)=-\xi(x{\cdot}v_j)=
-\sum_{i=1}^{n} c_{ij}\alpha_i.
$$
Since $x{\cdot}\xi=0$ if and only if $x{\cdot}\xi(v_j)=0$ 
for all $1\le j\le s$, 
the stabilisers $\gt q_\xi$ consists of all $x=\sum_{i}\alpha_ix_i$
such that 
$\sum_{i=1}^{n} c_{ij}\alpha_i=0$ for  $1\le j\le s$.
Hence $\dim\gt q_\xi=\dim\gt q-\rank A(\xi)$, 
where $A(\xi)$ is an $n{\times}s$-matrix with entries  $c_{ij}$.
Since $\rank A(\xi)\le\rank(\gt q{\cdot}V)$ and the equality holds for
generic $\xi\in V^*$, we get 
$\ind(\gt q,V)=\dim V^*-(\dim\gt q - (\dim\gt q- \rank(\gt q{\cdot}V))=
\dim V-\rank(\gt q{\cdot}V)$.
\end{proof}

In case $\gt q=\gt g_{0,e}$, $V=\gt g_{-1,e}$, we will denote 
the matrix $(\gt q{\cdot}V)$ by $([\gt g_{0,e},\gt g_{-1,e}])$.

\vskip0.5ex

Lemma~\ref{matrix} provides an easy method to 
compute an upper bound for the index that is very likely
to be equal to the index. 
Each $s$-tuple $\underline{a}=(a_1,\ldots,a_s)\in\ff^s$ defines an 
element 
%%% Let 
$\xi\in V^*$ such that  
%%  be fixed with 
$\xi(v_k) = a_k\in\ff$. 
Set $A(\underline{a}):=A(\xi)=(\xi(x_i{\cdot}v_j))$. 
Then  %%  we get
$$ \ind(\mf{q},V) = \dim V - \max_{\underline{a}\in\ff^s}
\rank( A(\underline{a})).$$
The entries of $A(\underline{a})$
are linear polynomials in the $a_i$. It follows that if we
take random coefficients $a_i$ then the rank of this matrix is very
likely maximal. 
%% We conclude that 
In other words,  
for  
%%% random choice of the 
any $s$-tuple $\underline{a}$ 
%% $a_k$ 
the
value of $\dim V -\rank A(\underline{a})$ is an upper bound for
$\ind(\mf{q},V)$, and if the $a_k$ are chosen randomly,
uniformly, and independently from a large enough set, then, very probably, equal to it.

There are several ways to get the value of the generic rank of 
$A(\underline{a})$. First of all, we can
consider the row space of $A(\underline{a})$ over the ground field
$\ff$, where we consider the $a_i$ as linearly independent indeterminates.
We can replace the rows by an $\ff$-linearly independent set of rows 
that span the same space over $\ff$. We can do the same with the columns.
Denote the resulting matrix by $\widetilde{A}(\underline{a})$. Then the
generic ranks of $A(\underline{a})$ and $\widetilde{A}(\underline{a})$ are the
same. If the lower bound that we get for the rank by substituting a point
$\underline{a}$ is equal to the number of columns, or rows, of 
$\widetilde{A}(\underline{a})$, then we know that this lower bound is 
the correct value of the generic rank. Otherwise we can compute the
rank of $\widetilde{A}(\underline{a})$, where the $a_i$ are indeterminates
of a function field over $\ff$. We do remark, however, that this operation
can be computationally expensive.

On the basis of Proposition~\ref{prop:2} we formulate the following algorithm.

\begin{alg}\label{alg_gib}
Input: a nilpotent elements $e\in\g_1$ and $\rank(G_0,\g_1)$.\\
Output: {\sc true} if $\rank(G_0,\g_1) = \ind(\g_{0,e},\g_{-1,e})$,
{\sc false} otherwise.

\begin{enumerate}
\item By linear algebra we compute bases of $\g_{0,e}$ and $\g_{-1,e}$. 
\item We compute the matrix $A(\underline{a})$ corresponding to the 
$\g_{0,e}$-module $\g_{-1,e}$.
\item By trying random values for the $a_i$ find a lower bound $r$ for the
generic rank of $A(\underline{a})$.
\item If $\dim \g_{-1,e}-r = \rank(G_0,\g_1)$ then output {\sc true}, 
else execute the next step. 
\item Compute $\widetilde{A}(\underline{a})$; if $r$ is equal to the
number of columns, or rows, of this matrix, then output {\sc false}, else
execute the next step.
\item Let $r'$ be the rank of $\widetilde{A}(\underline{a})$, where the
$a_i$ are indeterminates of a function field over $\ff$; 
if $\dim \g_{-1,e}-r' = \rank(G_0,\g_1)$ then output {\sc true}, else
output {\sc false}. 
\end{enumerate}
\end{alg}

\begin{lm}\label{lem:alggib}
The previous algorithm is correct.
\end{lm}

\begin{proof}
We note that $\dim \g_{-1,e}-r$ is an upper bound for $\ind(\g_{0,e},\g_{-1,e})$.
By Proposition \ref{prop:2}, $\rank(G_0,\g_1)$ is a lower bound for this
index. So if they are equal, then we know that we have the correct
index. If these are not equal, and $r$ is equal to the number of columns,
or rows, of $\widetilde{A}(\underline{a})$, then the index is strictly
bigger than $\rank(G_0,\g_1)$. Finally, if this also does not hold,
then the last ``brute force'' step gives the correct value.
\end{proof}

In order to check GIB for $(G_0,\g_1)$, we can do the following.
First we compute representatives of the nilpotent $G_0$-orbits in $\g_1$,
using the algorithms of \cite{gra15}. 
Since $\rank(G_0,\gt g_1)=\dim\gt g_1-\max_{e\in\gt g_1}\dim(G_0{\cdot}e)$, 
where $e$ is a nilpotent element, 
these calculations also provide the value of  $\rank(G_0,\gt g_1)$.
Then  
for each representative $e$ of a nilpotent $G_0$-orbit in $\g_1$
we execute Algorithm \ref{alg_gib}.
If the output is {\sc true} for all representatives of 
the nilpotent $G_0$-orbits in $\g_1$, then $(G_0,\g_1)$
has GIB by Lemma \ref{lem:alggib} and Proposition \ref{prop:2}.
If {\sc false} is returned once, then $(G_0,\g_1)$ does not have GIB.

\begin{rmk}
In practice it is a good idea to delay the execution of the expensive Step 6
of Algorithm \ref{alg_gib}. First one collects all nilpotent orbits for
which the random procedure indicates that GIB fails. Then Step 6 is executed
only once, on the smallest matrix $\widetilde{A}(\underline{a})$.
If the corresponding orbit does have GIB, which is very unlikely,
or the matrix is too complicated for computer calculations, 
one may look on other suspicious orbits. 
\end{rmk}

%%%%%%%%%%%%%%%%%%%%%%%%%%%%%%%%%%%%%%%%%%%
%%%%%%% SECTION %%%%%%%%%%%%%%%%%%%%%%%%%%% 
\section{Basic facts concerning semisimple inner automorphisms in type A}\label{basic-A}
%%%%%%%%%%%%%%%%%%%%%%%%%%%%%%%%%%%%%%%%%%%
%%%%%%%%%%%%%%%%%%%%%%%%%%%%%%%%%%%%%%%%%%%

Let $\VV$ be a finite dimensional vector space over $\ff$, 
$\gt g=\gt{gl}(\VV)$, and  
$\theta$ an inner automorphism of $\gt g$ of order $m$.
We consider $\theta$ as an element of the group $G=\GL(\VV)$ acting 
on $\gt g$ by conjugation.  Let $\zeta$ be a primitive $m$-th root of unity.
Set $\VV_t:=\{v\in \VV\mid \theta(v)=\zeta^t v\}$ and $r_t=\dim\VV_t$.
Up to a $G$-conjugation, $\theta$ is uniquely defined by the 
multiplicities vector $\hat r:=(r_0,r_1,\ldots,r_{m-1})$.  
Cyclic permutations of the entries of $\hat r$ correspond to 
multiplications by central elements of $\GL(\VV)$ and the resulting 
vectors $\hat r '$ define the same automorphism.  
We have $G_0=\GL(\VV_0){\times}\ldots{\times}\GL(\VV_{m-1})$ and 
$\gt g_1\cong\bigoplus\limits_{i=0}^{m-1}\Hom(\VV_i,\VV_{i+1})$, 
where $i+1$ is considered modulo $m$. 
Having $\hat r$, it is possible to write a Kac diagram of the corresponding 
$\theta$ and vice versa. Since we are not going to use this correspondence, 
it is only illustrated on one example. 

\begin{ex} Let $\theta$ be an automorphism of $\gt {gl}_9$
with the Kac diagram:

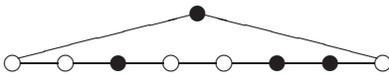
\begin{picture}(142,30)
  \put(3,0){\circle{6}}
  \put(23,0){\circle{6}}
  \put(43,0){\circle*{6}}
  \put(63,0){\circle{6}}
  \put(83,0){\circle{6}}
  \put(103,0){\circle*{6}}
  \put(123,0){\circle*{6}}
  \put(143,0){\circle{6}}

\put(73,19){\circle*{6}}
\put(5,2){\line(4,1){67}}
\put(141,2){\line(-4,1){67}}

  \put(6,0){\line(1,0){14}}
  \put(26,0){\line(1,0){14}}
  \put(46,0){\line(1,0){14}}
  \put(66,0){\line(1,0){14}}
%%  \put(43,3){\line(0,1){14}}
\put(86,0){\line(1,0){14}}
\put(106,0){\line(1,0){14}}
\put(126,0){\line(1,0){14}}
\end{picture} 

\vskip1ex

\noindent
Then $\theta$ is defined by $\hat r=(3,3,1,2)$.
\end{ex}

\begin{prop}(\cite[\S~3]{vinberg}) 
If $\theta$ is defined by $\hat r$, then 
$\rank(G_0,\gt g_1)=\min\limits_{i=0,\ldots,m-1} r_i$.
\end{prop}

Suppose that $x\in\gt g_1$. 
Let $x=x_s+x_n$ be the Jordan decomposition of $x$ in $\gt g$.
Due to its uniqueness, we have $x_s,x_n\in\gt g_1$. 
In other words, $\gt g_1$ inherits the 
Jordan decomposition from $\gt g$. 
The r{\^o}le of semisimple and nilpotent elements  in 
checking GIB is explained in  \cite[Section~2]{py-gib}. 
Suppose that  
$s\in\gt g_1$ is a semisimple element.
Then the action of $G_{0,s}$ on 
$\gt g_1/[\gt g_0,s]$ is called 
a {\it slice} representation of $(G_0,\gt g_1)$.
By the same argument as 
in the proof of Proposition~\ref{prop:2},
$(\gt g_1/[\gt g_0,s])^*\cong\gt g_{-1,s}$. 

\begin{lm}\label{slice}
Suppose that $\theta$ is defined by a vector 
$\hat r=(r_0,r_1,\ldots,r_{m-1})$ and 
the corresponding representation of $G_0$ on $\gt g_1$ has GIB. 
Let $b\in\mathbb N$ be such that $b\le\rank(G_0,\gt g_1)$.
Then the $\theta$-representation corresponding 
to $\hat r'=(r_0-b,r_1-b,\ldots,r_{m-1}-b)$ has GIB. 
\end{lm}
\begin{proof}
By inductive reasons, it suffices to prove the lemma for $b=1$.
Thus assume that $\rank(G_0,\gt g_1)>0$ and $b=1$.
Let $s\in\gt g_1$ be a non-zero semisimple element. %% such that 
%%% $\dim(s{\cdot}V_i)=1$ for all $i$. 
Set $r_i':=r_i-1$.
Then the slice representation of $G_{0,s}$ on $\gt g_{-1,s}^*$
is equivalent to the representation 
of 
$$
\ff^{^\times}\times\GL_{r_0'}{\times}\GL_{r_1'}{\times}\ldots{\times}\GL_{r_{m-1}'}
\ \text{ on } \
\ff\oplus\left(\bigoplus\limits_{i=0}^{m-1}(\ff^{{r_i}'})^*{\otimes}\ff^{ {r_{i+1}}'}\right),
$$
where the first subgroup $\ff^{^\times}$ acts on $\gt g_{-1,s}^*$
trivially and the first subspace $\ff$ is a trivial $G_{0,s}$-module.  
It follows that $(G_{0,s},\gt g_{-1,s}^*)$ has GIB if and only if the 
$\theta$-representation corresponding to $\hat r'$ has GIB. It remains to 
notice that GIB of $(G_0,\gt g_1)$ implies GNIB 
of $(G_{0,s},\gt g_{-1,s}^*)$ by \cite[Theorem~2.1]{py-gib}
and in our setting GNIB is equivalent to GIB by \cite[Theorem~2.3]{py-gib}.
\end{proof}

\subsection{Nilpotent orbits.}\label{nilp-g1}
Nilpotent elements in $\gt{gl}(\VV)$ are classified in terms of their Jordan 
normal form. A similar classification is possible for nilpotent 
$G_0$-orbits in $\gt g_1$. This problem was solved by Kempken in 
\cite[\S 2.II]{Gisela-dis}. 
Below we present her answer in a slightly modified form. 

Let $e\in\gt g_1$ be a nilpotent element represented by a 
partition $(d_1{+}1,\ldots,d_k{+}1)$ of $\dim\VV$.
Since the subspace $e{\cdot}\VV\subset\VV$ is $\theta$-invariant 
and $\theta$ acts on $\VV$ as  
a semisimple element, there is  
a $\theta$-invariant complement to $e{\cdot}\VV$, let us say, $\W$.  
Clearly $\dim\W=k$.
Assume that $d_1\ge d_2\ge\ldots\ge d_k$.  
For each number $s$ the subspaces 
$\ker e^s=\{v\in\VV\mid e^s{\cdot}v=0\}$ and $\ker e^s\cap \W$ are both 
$\theta$-invariant. Hence 
a generator $w_1$ of the maximal Jordan block, i.e.,
a vector such that $e^{d_1}{\cdot}w_1\ne 0$, can be chosen
as an eigenvector of $\theta$. Set 
$\VV[1]:=\textrm{span}\{w_1,e{\cdot}w_1,\ldots,e^{d_1}{\cdot}w_1\}$.
Let $\W'$ be a $\theta$-invariant complement of $\ff w_1$ in 
$\W$. Then also $\VV=\VV[1]\oplus (\W'{\oplus}e{\cdot}\W')$ is 
a $\theta$-invariant decomposition. 
Proceeding by induction on the number of Jordan 
blocks we prove that all generators $w_i$ can be chosen as
eigenvectors of $\theta$, i.e., $\W$ has a basis 
$w_1, w_2, \ldots, w_k$ consisting of  
$\theta$-eigenvectors, where in addition   
the vectors $e^{j}{\cdot}w_i$ with 
$1\le i\le k$, $0\le j\le d_i$ form a basis of $\VV$.  
Note that $e^{d_i+1}{\cdot}w_i=0$ for all $i\le k$. 
If $\theta(w_i)=\zeta^{t(i)}w_i$, then
$\theta(e^s{\cdot}w_i)=\theta(e^s){\cdot}\theta(w_i)= 
 \zeta^{t(i)+s}e^s{\cdot}w_i$. 
In particular,  $\theta(e{\cdot}w_i)=\zeta^{t(i)+1}e{\cdot}w_i$.           
Summing up:
\begin{itemize}
\item[ ] \ 
to each nilpotent element $e\in\gt g_1$ we associate its partition and 
the $\theta$-eigenvalues $\zeta^{t(i)}$ on the Jordan blocks' generators 
$w_i$. 
\end{itemize}
In order to see what nilpotent elements do appear in 
$\gt g_1$, we have to take a partition, choose $\theta$-eigenvalues
for the generators $w_i$ and count the dimensions of
the eigenspaces according to the rule 
$\theta(e^s{\cdot}w_i)=\zeta^{t(i)+s}e^s{\cdot}w_i$. 
If they coincide with the $r_t$'s, then $e$ lies in $\gt g_1$.

\subsection{Basis of a centraliser.}\label{basis}
In order to do explicit calculations,
one needs bases in $\gt g_{0,e}$ and $\gt g_{-1,e}$.
First we introduce a basis in $\gt g_e$. 
If $\xi\in\g_e$, then $\xi{\cdot}(e^j{\cdot}w_i)=e^j{\cdot}(\xi{\cdot}w_i)$, hence
$\xi$ is completely determined by its values on $\W$. 
The only restriction on $\xi{\cdot}w_i$ is that 
$e^{d_i+1}{\cdot}(\xi{\cdot}w_i)=\xi{\cdot}(e^{d_i+1}{\cdot}w_i)=0$. 
Since vectors $e^s{\cdot}w_i$ form a basis of $\VV$, the
centraliser $\gt g_e$ has a basis 
$\{\xi_i^{j,s}\}$ such that
$$
\left\{
\begin{array}{l}
\xi_i^{j,s}{\cdot}w_i=e^s{\cdot}w_j, \\
\xi_i^{j,s}{\cdot}w_t=0 \enskip \mbox{for } t\ne i, \\
\end{array}\right.
\quad 1\le i,j\le k, \ \mbox{ and }\ \max\{d_j-d_i, 0\} \le s\le d_j \ .
$$
It is convenient to assume that $\xi_i^{j,s}=0$ whenever 
$s$ does not satisfy the above restrictions. 
The composition rule shows that the basis elements 
$\xi_i^{j,s}$ satisfy the following commutator relation:
\begin{eqnarray}\label{commutator}
[\xi_i^{j,s},\xi_p^{q,t}]=\delta_{q,i}\xi_p^{j,t+s}-\delta_{j,p}\xi_i^{q,s+t},
\end{eqnarray}
where $\delta_{i,j}=1$ if $i=j$ and is zero otherwise. 
%%see \cite{fan} for more detail. 
Each $\xi_i^{j,s}$ is an eigenvector of $\theta$ with 
%% More pricesely, 
\begin{equation}\label{sobst-zn}
\theta(\xi_i^{j,s})=\zeta^{s+t(j)-t(i)}\xi_i^{j,s}.
\end{equation}
This allows one to compute $\gt g_{0,e}$ and $\gt g_{-1,e}$.

\begin{ex}\label{nilp-ex}
Let $\theta$ be an automorphism of $\gt{gl}_9$ 
with $\hat r=(3,3,3)$. Then there is a nilpotent 
element $e\in\gt g_1$ defined by a partition 
$(5,3,1)$ such that $\theta(w_1)=w_1$, 
$\theta(w_2)=\zeta w_2$, and $\theta(w_3)=\zeta^2 w_3$.
Indeed, let us put the eigenvalues of $\theta$, or rather 
the corresponding exponents of $\zeta$, in the squares of the
Young diagram corresponding to $e$.
\begin{center}
\begin{picture}(45,80)(0,0)

\put(0,0){\line(1,0){45}}
\put(0,15){\line(1,0){45}}
\put(0,30){\line(1,0){30}}
\put(0,45){\line(1,0){30}}
\put(0,60){\line(1,0){15}}
\put(0,75){\line(1,0){15}}

\put(0,0){\line(0,1){75}}
\put(15,0){\line(0,1){75}}
\put(30,0){\line(0,1){45}}
\put(45,0){\line(0,1){15}}

\put(5,4){$0$}
\put(5,19){$1$}
\put(5,34){$2$}
\put(5,49){$0$}
\put(5,64){$1$}

\put(20,4){$1$}
\put(20,19){$2$}
\put(20,34){$0$}

\put(35,4){$2$}
\end{picture}
\end{center}
Then one readily sees that there are three $0$, three $1$, 
and three $2$ in the figure.  
The centraliser $\gt g_e$ has a basis
$$
\xi_1^{1,0},\xi_1^{1,2},\xi_1^{1,2},\xi_1^{1,3},\xi_1^{1,4}, \ \ \xi_1^{2,0},\xi_1^{2,1},\xi_1^{2,2}, \ \
\xi_1^{3,0}, \ \ \xi_2^{2,0},\xi_2^{2,1},\xi_2^{2,2}, \ \ \xi_2^{1,2},\xi_2^{1,3},\xi_2^{1,4}, \ \ \xi_2^{3,0}, \ \
\xi_3^{3,0},\xi_3^{1,4},\xi_3^{2,2},
$$
where the $\theta$-eigenvalues are
$$
1,\zeta,\zeta^2,1,\zeta,\ \ \zeta,\zeta^2,1,\ \ \zeta^2,\ 1,\zeta,\zeta^2,\ \  \zeta,\zeta^2,1, \ \  \
 \zeta,1,\zeta^2,\zeta,
$$
respectively. In particular $\dim\gt g_{0,e}=6$, $\dim\gt g_{1,e}=7$, and
$\dim\gt g_{-1,e}=6$. 
\end{ex}
 
\section{GIB in type A}

In this section we consider inner finite order automorphisms 
$\theta$ of $\gt{gl}_n$. 
All $\theta$ such that $\rank(G_0,\gt g_1)>0$ and 
the pair $(G_0,\gt g_1)$ has GIB are classified. 
According to \cite{py-gib}, there are 
only three such involutions, 
namely $(\gt g,\gt g_0)$ must be one of the following pairs: 
$(\gt{gl}_{n+2},\gt{gl}_n{\oplus}\gt{gl}_2)$, 
$(\gt{gl}_{n+1},\gt{gl}_n{\oplus}\gt{gl}_1)$,
$(\gt{gl}_6,\gt{gl}_3{\oplus}\gt{gl}_3)$.  
Not surprisingly, 
for automorphisms of higher orders the GIB property can be satisfied only if 
$\rank(G_0,\gt g_1)\le 2$. 
All initial, so to say, nilpotent orbits without GIB 
were found on computer. After that it is possible to 
extend these bad examples to higher dimensions. 
We also check on computer that some 
$\theta$-representations of small dimension do have GIB. 
The computer calculations were done using our implementation of 
Algorithm~\ref{alg_gib}.

The difference between $\theta$-representations of 
$\gt{sl}_n$ and $\gt{gl}_n$ is almost neglectable.
Sometimes the general linear algebra is more convenient for calculations. 
On the other hand, we always have to take into account the central torus, 
which acts on $\gt g$ and $\gt g_1$ trivially. From now on let 
$z$ be a central element in $\gt{gl}_n$. 
We deal with automorphisms $\theta$ in terms of 
the corresponding vectors $\hat r$, as defined in Section~\ref{basic-A}.

Although the goal is to classify automorphisms of positive rank having GIB,
we first give an example of a $\theta$-representation with 
$\rank(G_0,\gt g_1)=0$. 

\begin{ex}\label{3,0} Suppose that $m=3$ and $\theta$ has rank zero, 
i.e., $\hat r=(a,b,0)$ up to a cyclic permutation. 
Then $(G_0,\gt g_1)$ has GIB. Indeed, 
$\gt  g_1=M_{a,b}(\mathbb F)$ and $G_0$-orbits $G_0{\cdot}x\subset\gt g_1$
are classified by the rank  $p$ of an $a{\times}b$-matrix $x$. The quotient 
space
$\gt g_1/[\gt g_0,x]$ is isomorphic to $M_{a-p,b-p}(\mathbb F)$ and 
$G_{0,x}$ acts on it as $\GL_{a-p}{\times}\GL_{b-p}$, hence, with an open orbit. 
%%%% and also with an open orbit in the dual space.  
\end{ex}

The $\theta$-representation of Example~\ref{3,0}
appears as a slice representation 
for the action of $G_0$ on $\gt g_1$ corresponding 
to $\hat r=(a+1,b+1,1)$. This second $\theta$-representation 
has GIB as well. 
In order to prove it, we need 
three following lemmas. 
%%% an auxilary statement in rank zero. 

\begin{lm}\label{gl_q}
Suppose that $\hat r=(1,n,1,0)$. Then the 
corresponding representation of $G_0$ on $\gt g_1$
has GIB. 
\end{lm}
\begin{proof}
Here the action of $G_0$ on $\gt g_1=\ff^n\oplus(\ff^n)^*$ 
has a one-dimensional 
ineffective kernel, say $Q_0$, and 
$H:=G_0/Q_0=\ff^{^\times}{\times}\SL_n{\times}\ff^{^\times}$.
We have $(t_1,t_2){\cdot}(v_1+v_2)=t_1v_1+t_2v_2$ for 
$(t_1,t_2)\in \ff^{^\times}\!{\times}\,\ff^{^\times}$,
$v_1\in\ff^n$, $v_2\in(\ff^n)^*$. Set $\gt h=\Lie H$.
According to Proposition~\ref{summa}, 
we have to show that for all $v=v_1+v_2\in\gt g_1$, there is 
$w\in\gt g_1$ such that 
$\gt h_v{\cdot}w+\gt h{\cdot}v=\gt g_1$. 
In cases $v=0$, where $\gt h_v{\cdot}w=\gt g_1$ for generic $w\in\gt g_1$;
and $v_2(v_1)\ne 0$, where $\gt h{\cdot}v=\gt g_1$,
the statement is clear. 
If one of the vectors $v_1,v_2$ is zero, 
without loss of generality we may assume that $v_2=0$, 
and the other one is not (now $v_1\ne 0$), then 
$\gt h{\cdot}v=\ff^n$ and $\gt h_v{\cdot}w_2=(\ff^n)^*$ 
for generic $w_2\in (\ff^n)^*$. 
 
It remains to treat the case where $v_1,v_2$ are both non-zero, 
but $v_2(v_1)=0$. This implies that $n\ge 2$.
Here $\dim(\gt h{\cdot}v)=2n-1$ and 
if $w_2\in(\ff^n)$ is such that $w_2(v_1)$, then 
$\ff w_2\cap \gt h{\cdot}v=\{0\}$. Let 
$\rho(\ff^{^\times})\in\SL_n$ be
a one-dimensional torus such that 
$\rho(t){\cdot}(v_1+v_2)=tv_1+tv_2$. 
Take an element $g_t=(t^{-1},\rho(t),t^{-1})\in H$.
Then $g_t{\cdot}w_2\in t^{-2}w_2+\gt h{\cdot}v$
and therefore $\gt h_v{\cdot}w_2$ contains $\ff w_2$.
We conclude that $\gt h_v{\cdot}w_2+\gt h{\cdot}v=\gt g_1$.
\end{proof}

\begin{lm}\label{rank0-2gr}
Suppose that $\hat r=(1,a,b,1,0)$. Then the 
corresponding representation of $G_0$ on $\gt g_1$
has GIB. 
\end{lm}
\begin{proof}
Using explicit matrix calculations
we check that all orbits in $\gt g_1$ satisfy the 
GIB property. Here $G_0=\ff^{^\times}{\times}\GL_a{\times}\GL_b{\times}\ff^{^\times}$ and
there is a $G_0$-invariant decomposition 
$\gt g_1=V_1\oplus V_2\oplus V_3$, where   
$V_1\cong (\ff^a)^*$, $V_2\cong \ff^a{\otimes}(\ff^b)^*\cong M_{a,b}(\ff)$, 
and $V_3\cong \ff^b$. 
Take $e\in\gt g_1$ and let $x\in M_{a,b}(\ff)$ be its
projection on $V_2$ (this means that $e\in x+V_1+V_3$).  
Let $q$ be the rank of the matrix $x$. 
Replacing $e$ by another element 
in  $G_0{\cdot}e$ we may (and will) assume that $x$ is an identity 
$q{\times}q$ matrix standing in the upper left corner. 
Then  the stabiliser $G_{0,x}$ is a product 
$\ff^{^\times}{\times}\GL_{a-q}{\times} U_a{\times}\GL_q{\times} U_b{\times} \GL_{b-q}{\times}\ff^{^\times}$,
where $\GL_q$ is embedded diagonally into $\GL_a{\times}\GL_b$ and
$U_a,U_b$ are unipotent radicals of standard parabolics 
in $\GL_a,\GL_b$, respectively.  Set $\gt u_a=\Lie U_a$, 
$\gt u_b=\Lie U_b$. 

Now $V_1=(\ff^q)^*\oplus(\ff^{a-q})^*$, 
$V_3=\ff^q\oplus\ff^{b-q}$, where 
$\GL_q$ acts non-trivially only on $\ff^q\oplus(\ff^q)^*$, 
the subgroup $\GL_{a-q}$ only on $(\ff^{a-q})^*$,
and $\GL_{b-q}$ only on $\ff^{b-q}$.
For the nilpotent radicals we have  
$[\gt u_a,(\ff^q)^*]=(\ff^{a-q})^*$, if $q\ne 0$, while 
$[\gt u_a,(\ff^{a-q})^*]=0$ and likewise
$[\gt u_b,\ff^q]=\ff^{b-q}$, if $q\ne 0$, while $[\gt u_b,\ff^{b-q}]=0$. 
According to this decomposition, we write $e$ is a sum 
of five vectors $e=v_1+v_1'+x+v_3+v_3'$ with 
$v_1\in (\ff^{q})^*$, $v_3\in\ff^q$.
Set $W=\ff^q\oplus(\ff^q)^*$.
The reductive part 
of $G_{0,x}$ acts on $W$ in exactly the same way
as the $\theta$-group in Lemma~\ref{gl_q}.
Let $H$ be its image in $\GL(W)$ and $\gt h:=\Lie H$.   
Then for some vector $w\in W$, and therefore 
for all vectors of an open subset, holds 
$\gt h_{v_1{+}v_3}{\cdot}w+\gt h{\cdot}(v_1{+}v_3)=W$. 
There are three different possibilities, which are treated separately. 

\noindent
{\bf 1.} \ Suppose that $v_1$ and $v_3$ are both non-zero. 
Replacing, if necessary, $e$ by an element in $(U_a{\times}U_b){\cdot}e$,
we may assume that $v_1'=v_3'=0$.
%%% In case $q>1$, for generic $w=w_1+w_3\in W$ we have 
%%% $[(\gt u_a)_{v_1},w_1]=(\ff^{a-q})^*$ and 
%%% $[(\gt u_b)_{v_3},w_3]= \ff^{b-q}$. 
Let $w\in W$ be a generic vector and 
$y\in V_2$ a matrix such that its 
lower right $(a-q){\times}(b-q)$ submatrix  
is of the full rank $\min(a-q,b-q)$. 
Then  
%%$$
%%\begin{array}{l}
%%
\begin{equation}\label{1.}  
 \begin{split}
 & [\gt g_{0,e},w+y]+[\gt g_0,e]=[\gt g_{0,e},w+y]+[\gt g_{0,x},v_1+v_3]+[\gt g,e]=  \\
& \enskip %%%% [(\gt u_a)_{v_1},w_1]+[(\gt u_b)_{v_3},w_3]+
 [\gt{gl}_{a-q}{\oplus}\gt{gl}_{b-q},y]+ \gt h_{v_1+v_3}{\cdot}w +
\gt h{\cdot}(v_1{+}v_3)+[\gt u_a{\oplus}\gt u_b,v_1+v_3]+[\gt g_0,e]=  \\
 & =W+(\ff^{q-a})^*+\ff^{q-b} + [\gt{gl}_{a-q}{\oplus}\gt{gl}_{b-q},y]+[\gt g_0,x]=\gt g_1.   
\end{split}
\end{equation}
%%\end{array}
%%$$
%%% In case $q=1$, we have $[\gt g_0,e]=V_1+V_3+[\gt g_0,x]$.
%%% Hence $\gt g_1/[\gt g_0,e]\cong (\ff^{a-1})^*{\otimes}\ff^{b-1}$
%%% and $G_{0,e}=\GL_{a-1}{\times}\GL_{b-1}$ acts on it with an open 
%% orbit.  
By Proposition~\ref{summa}, the element $e$ has GIB. 

\noindent
{\bf 2.} \ Suppose now that one of the vectors $v_1,v_3$ is zero, 
but the other one is not. Without loss of generality we 
may assume that $v_1\ne0$, $v_3=0$. Now $v_3'$ cannot be assumed 
to be zero, but $\gt u_b\subset\gt g_{0,e}$. 
In the equation~(\ref{1.}) we have to replace
$\gt{gl}_{b-q}$ by $(\gt{gl}_{b-q})_{v_3'}$ in 
$[\gt{gl}_{a-q}{\oplus}\gt{gl}_{b-q},y]$;
and $[\gt u_b,v_3]$ by $[\gt u_b,w]$, which 
is equal to $\ff^{b-q}$ for generic $w$ and is a subset 
of $[\gt g_{0,e},w+y]$. 
Finally notice that still 
$[\gt{gl}_{a-q}{\oplus}(\gt{gl}_{b-q})_{v_3'},y]+[\gt g_0,x]=V_2$.
Therefore $[\gt g_{0,e},w+y]+[\gt g_0,e]=\gt g_1$.

\noindent
{\bf 3.} \ The last possibility is that $v_1=v_2=0$. 
If $v_1'\ne 0$ and 
$v_3'\ne 0$, then
$[(\gt{gl}_{a-q})_{v_1'}{\oplus}(\gt{gl}_{b-q})_{v_3'},y]+[\gt g_0,x]$
is a subspace of codimension $1$ in $V_2$.
Otherwise the sum is the whole of $V_2$ and again 
$\gt g_1=[\gt g_{0,e},w+y]+[\gt g_0,e]$.
Thus we may safely assume that both vectors are non-zero. 
In particular, $V_1{\oplus}V_3\subset[\gt g_0,e]+V_2$. 

We have $\gt g_1/[\gt g_0,e]=\gt w_1{\oplus}\gt w_2{\oplus}\gt w_3$,
where 
$\gt w_2\subset V_2$, $\gt w_2\cong M_{a-q,b-q}(\ff)$, 
$\gt w_1\cong (\ff^q)^*$ is embedded anti-diagonally into 
$V_1{\oplus} V_2$, and 
$\gt w_1\cong \ff^q$ is embedded anti-diagonally into 
$V_2{\oplus} V_3$. 
Next   
$$
\begin{array}{l}
G_{0,e}=\GL_q\times U_a\times U_b \times 
(\ff^{^\times}{\times}\GL_{a-q}{\times}\GL_{b-q}{\times}\ff^{^\times})_{(v_1'+v_2')}. \\
\end{array}
$$
Let $y\in\gt w_2$ be a matrix of the maximal rank. 
Then $\gt u_a{\cdot}y=\gt w_1$, $\gt u_b{\cdot}y=\gt w_3$, 
and $\gt g_{0,e}{\cdot}y=\gt g_1/[\gt g_0,e]$.
Thereby $\ind(\gt g_{0,e},(\gt g_1/[\gt g_0,e])^*)=0$. 
\end{proof}

\begin{lm}\label{fold}
Let  $\theta$ be an automorphism of a Lie algebra 
$\gt g=\gt{gl}(\VV)$ of order $m$
defined by  a vector $\hat r=(1,r_1,\ldots,r_{m-1})$.
Suppose that $\theta'$ is an automorphism of order $m+1$ of another 
Lie algebra $\gt h=\gt{gl}(\VV')$, where $\dim\VV'=\dim\VV+1$,
defined by a vector $\hat r'=(1,r_1,\ldots,r_{m-1},1,0)$,
and $\theta'$ has GIB.
Then  $\theta$ has GIB as well. 
\end{lm}
\begin{proof}
Note that $\gt g_1=\gt h_1$. The group $G_0$ is smaller than 
$H_0$ and does not contain the central torus of $H_0$ (acting on 
$\gt g_1$ by the scalar multiplications).
Therefore we are in the setting of Lemma~\ref{ind-0-1}
and the result follows  from it. 
\end{proof}

In case $m=3$ we can get a complete answer. 
This is achieved in a few following steps.  
We will need a machinery developed in subsections~\ref{nilp-g1} 
and \ref{basis}.

\begin{prop}\label{3-1}
Suppose that $m=3$ and $\theta$ has rank one,
i.e., $\hat r=(a,b,1)$  with $a,b>0$. %%% up to a cyclic permutation. 
Then $(G_0,\gt g_1)$ has GIB.
\end{prop}
\begin{proof}
We get the result combining Lemmas~\ref{rank0-2gr} and \ref{fold}. 
\end{proof}

\begin{lm}\label{r=2;4} Suppose that $\hat r=(2,2,a)$ with 
$0\le a\le 4$. Then the corresponding automorphism   
$\theta$ has GIB.
\end{lm}
\begin{proof}
For $a=0,1$ the statement follows from Example~\ref{3,0} and Proposition~\ref{3-1}.
In cases $a=2,3,4$, GIB was checked according to Algorithm~\ref{alg_gib}.    
Direct verification by hand is also possible. 
\end{proof}

\begin{prop}\label{m=3;2,a} Suppose that $\hat r=(2,2,a)$.
Then $\theta$ has GIB for all $a$.
\end{prop}
\begin{proof}
Due to Lemma~\ref{r=2;4}, the statement is true 
for $a\le 4$. Assume that $a>4$. 
Suppose that $\gt g=\gt{gl}(\VV)$ with $\dim\VV=a{+}4$. 
Let $\gt h=\gt{gl}(\ff^8)$ with $\ff^8\subset\VV$ 
be a $\theta$-invariant subalgebra 
of $\gt g$  such that 
the restriction of $\theta$ to $\gt h$ is an automorphism 
with $\hat r_{\gt h}=(2,2,4)$.
Let also $H\subset G$ be a connected subgroup with $\Lie H=\gt h$. 
Take a nilpotent element $e\in\gt g_1$. 
It must have at least $a-4$ Jordan blocks of size 
zero such that $\theta(w_i)=\zeta^2 w_i$.
Therefore $G_0{\cdot}e\cap \gt h_1\ne\varnothing$
and we may (and will) assume that $e\in\gt h_1$.  

%%% We need to understand 
%%% $\gt g_1/[\gt g_0,e]$ and its relation with 
%%  $\gt h_1/[\gt h_0,e]$. 
Let $\gt m$ stand for the $\gt h$-invariant complement 
of $\gt h_1$ in $\gt g_1$. As a linear space 
$\gt m=M_{a-4,2}(\ff)\oplus M_{2,a-4}(\ff)$.
Let $\gt p,\gt p_-\subset\gt{gl}_a\subset\gt g_0$ 
be two opposite parabolic
subalgebras with the Levy part $\gt{gl}_4\oplus\gt{gl}_{a-4}$
and $\gt u_1,\gt u_2$ their nilpotent radicals. 
Note that $[\gt u_1,\gt u_1]=[\gt u_2,\gt u_2]=0$
and $\dim\gt u_1=\dim\gt u_2=4(a-4)$.
Next 
$\gt h_1=V_1\oplus V_2\oplus V_3$, 
where $V_1=M_{4,2}(\ff)$, $V_2=M_{2,2}(\ff)$, 
and $V_3=M_{2,4}(\ff)$ are $H_0$-invariant subspaces.
Let us write $v=x_1+x_2+x_3$ in  accordance  with this decomposition 
of $\gt h_1$. 
We suppose that $\gt u_1$ is presented by 
above-diagonal matrices and $\gt u_2$ by below-diagonal. 
Then $[\gt u_1, x_1]=[\gt u_2,x_3]=0$ and $[(\gt u_1{\oplus}\gt u_2),x_2]=0$.  
Apart from this we have  
$[(\gt u_1)_{x_3},w]=M_{2,a-4}(\ff)$ for generic 
$w\in V_3$. 
Similar equality holds for $\gt u_2$. 
As a consequence, for generic $w\in\gt h_1$, the subspace 
$[(\gt u_1{\oplus}\gt u_2)_e,w]$ coincides with $\gt m$. 
Hence 
$$
[\gt g_{0,e},w]+[\gt g_0,e]\supset ( [\gt h_{0,e},w] + [\gt h_0,e] )\oplus 
 [(\gt u_1)_{e},w] \oplus [(\gt u_2)_e,w] =
  ( [\gt h_{0,e},w] + [\gt h_0,e] ) \oplus \gt m
$$ 
for a generic $w\in\gt h_1$. 
By Lemma~\ref{r=2;4}, the pair $(H_0,\gt h_1)$ has GIB. 
Taking the intersection of two open subsets of $\gt h_1$, 
we may assume that $w$ yields also an $H_{0,e}$-orbit of codimension $2$ 
in $\gt h_1/[\gt h_0,e]$. Then  
$G_{0,e}$-orbit of $w$ is of codimension at most $2$ in 
$\gt g_1/[\gt g_0,e]$. Hence the pair $(G_0,\gt g_1)$ has GIB as well. 
\end{proof}

\begin{ex}\label{ex-3,3,2} Let $\theta$ be an automorphism of 
$\gt{gl}_8$ of order $3$ with $\hat r=(3,3,2)$. 
Then there is no GIB here. Indeed, 
take a nilpotent element 
$e\in\gt g_1$ having Jordan blocks $(5,3)$ with 
$\theta(w_1)=w_1$ 
and $\theta(w_2)=\zeta w_2$. %%%%  labeled as $1$. 
We can choose bases of 
 $\gt g_{0,e}$ and  $\gt g_{-1,e}$ as follows:
%%%%is $4$-dimensional and 
%%%% has a basis
$$
z, \xi_1^{1,0},\xi_1^{1,3}, \xi_1^{2,2}, \xi_2^{1,4}; \qquad
\xi_1^{1,2},\xi_1^{2,1},\xi_2^{1,3},\xi_2^{2,2}.
%%% 
%% $\gt g$ anyway, $\gt g_{0,e}$ is also $4$-dimensional 
%% with a basis
$$
The nilpotent radical of $\gt g_{0,e}$ is three dimensional and is generated by the last three basis vectors. By (\ref{commutator}) it
commutes with $\gt g_{-1,e}$, 
for example, $[\xi_1^{2,2},\xi_2^{1,3}]=\xi_2^{2,5}-\xi_1^{1,5}=0$.
Therefore the matrix $([\gt g_{0,e},\gt g_{-1,e}])$ has rank at most $1$ and,
due to Lemma~\ref{matrix}, 
$\ind(\gt g_{0,e},\gt g_{-1,e}^*)\ge 3>2$. 
\end{ex}

In the next example,
among 191 nilpotent 
$G_0$-orbits in $\gt g_1$ 
there are three bad ones, without GIB. One of them is 
presented here. Others arise after cyclic permutations of 
$\theta$-eigenvalues on $w_1,w_2,w_3$. 

\begin{ex}\label{ex-3,3,3} Let $\theta$ be an automorphism of 
$\gt{gl}_9$ of order $3$ with $\hat r=(3,3,3)$. 
Consider a nilpotent element
$e\in\gt g_1$ with Jordan blocks of sizes $(5,3,1)$, where 
$w_1$ is $\theta$-invariant, $\theta(w_2)=\zeta w_2$, and 
$\theta(w_3)=\zeta^2 w_3$ (the same as in Example~\ref{nilp-ex}).
%%%%%%%% defined by the labeled partition $([5,0],[3,1],[1,2])$.  
%%% Then $\gt g_{0,e}$ is generated by three 
%%% semisimple elements (one of which is central) and three other vectors 
The subspaces $\gt g_{0,e}$ and $\gt g_{-1,e}$ have bases 
$$
z,\xi_2^{2,0},\xi_3^{3,0},\xi_1^{1,3}, \xi_1^{2,2}, \xi_2^{1,4}; 
\ \text{ and } \ \ 
\xi_1^{1,2},\xi_1^{2,1},\xi_1^{3,0},\xi_2^{1,3},\xi_2^{2,2},\xi_3^{1,4},
$$ 
respectively. 
Using (\ref{commutator}) it is not difficult to verify  that the 
nilpotent part of $\gt g_{0,e}$, generated by $\xi_1^{1,3}, \xi_1^{2,2}, \xi_2^{1,4}$,  
commutes with $\gt g_{-1,e}$.
Hence $\ind(\gt g_{0,e},\gt g_{-1,e})\ge 6-2=4>3=\ind(\gt g_0,\gt g_1)$.
\end{ex} 

These examples lead to the following statement. 

\begin{prop}\label{m=3-prop}
Suppose that $\hat r=(r_0,r_1,r_2)$ with 
$r_0,r_1>2$ and $r_2\ge 2$. Then $(G_0,\gt g_1)$ does not 
have GIB. 
\end{prop}
\begin{proof}
Due to Lemma~\ref{slice},
it suffices to prove the statement for $\hat r$ such that 
$\hat r'=(r_0-1,r_1-1,r_2-1)$ does not satisfy the assumptions after 
any cyclic permutation. This is possible an exactly two cases:
$r_2=2$ or if at least two of the numbers $r_0,r_1,r_2$ are equal to $3$.  
%%
%%
%%
%%%%%%% 
 
Suppose first that $\hat r=(a,b,2)$ with $a,b>2$. 
Let $e\in\gt g_1$ be the same nilpotent element as
in Examples~\ref{ex-3,3,2},\ref{ex-3,3,3}, i.e., 
$e$ has Jordan blocks of sizes $(5,3,1^{a+b-6})$. 
For the generators of these Jordan blocks holds:
$\theta(w_1)=w_1$, $\theta(w_2)=\zeta w_2$,
$\theta(w_i)=w_i$ for $3\le i\le a-1$, 
and $\theta(w_j)=\zeta w_j$ for $a\le j\le a+b-2$. 
Let $\ff^8\subset\VV$ be a $\theta$-invariant subspace 
such that the restriction of $\theta$ to $\gt h:=\gt{gl}(\ff^8)$
is defined by $\hat r_{\gt h}=(3,3,2)$. 
Let $\gt f\cong\gt{gl}_{a{+}b{-}6}$ be a $\theta$-invariant 
subalgebra such that  
$\gt h{\oplus}\gt f\subset\gt g$ is a 
Levi subalgebra of $\gt g$.  
Set $\gt a:=\gt h_{-1,e}=\gt h\cap\gt g_{-1,e}$.
Then $\dim\gt a=4$ and as a vector space $\gt a$ is 
generated by 
$\xi_1^{1,2},\xi_1^{2,1},\xi_2^{1,3}$, and $\xi_2^{2,2}$.
Let $\gt m\subset\gt g_{0,e}$  be a subspace generated
by $\xi_i^{j,s}\in\gt g_{0,e}$ such that either $i=1,2$ and $j>2$ (then 
necessary $s=0$), or $i>2$ and $j=1,2$
(then necessary $s=d_j$). 
Then $\gt g_{0,e}=\gt h_{0,e}{\oplus}\gt m{\oplus}\gt f_0$.
%%%% where $f\cong\gt{gl}_{a{+}b{-}6}$ is a commuting with $\gt h$ 
%%%% subalgebra of $ $. 
Suppose $\eta\in\gt m$. 
Among 
$\xi_1^{1,t},\xi_1^{2,t},\xi_2^{1,t},\xi_2^{2,t}\in\gt h_e$ 
the element $\eta$ does not commute only with 
$\xi_1^{1,0},\xi_1^{2,0},\xi_2^{1,2},\xi_2^{2,0}$. 
In particular, $[\gt a,\eta]=0$. Note also that $[\gt a,\gt f]=0$. 
Hence
$\rank([\gt g_{0,e},\gt a])=\rank([\gt h_{0,e},\gt a])\le 1$
(the  inequality  was shown to be true in Example~\ref{ex-3,3,2}). 
%%
%%This means that %%%% the first four lines of the matrix 
%%% $([\gt g_{-1,e},\gt g_{0,e}])$, corresponding to 
%%%%% $([\gt a,\gt g_{0,e}])$, are proportional to each other 
Due to Lemma~\ref{matrix},
$$
\ind([\gt g_{0,e},\gt g_{-,e}])= \dim\gt g_{-1,e}-\rank([\gt g_{0,e},\gt g_{-1.e}])\ge 
\dim\gt a-\rank([\gt g_{0,e},\gt a])\ge 3.
$$ 
Therefore $\ind(\gt g_{0,e},\gt g_{-1,e})>\ind(\gt g_0,\gt g_1)=2$. 

Now we pass to the second case, where 
$r_2\ge 3$ and at least two of the 
numbers $r_0,r_1,r_2$ are equal to $3$.
Without loss of generality we may assume that
$\hat r=(3,3,a)$ with $a\ge 3$. 
Here we use the same 
subalgebra $\gt h\subset\gt g$ and 
almost the same nilpotent element $e\in\gt h_1$.
The only difference is that now 
$e$ has $a{-}2$ Jordan blocks of size $1$ and for them holds
$\theta(w_i)=\zeta^2 w_i$, if 
$3\le i\le a$.
We have 
$\gt g_{0,e}=\gt h_{0,e}{\oplus}\gt{gl}_{a-2}$ and
$\gt g_{-1,e}=\gt h_{-1,e}\oplus\ff^{a-2}\oplus(\ff^{a-2})^*$, 
where $\gt{gl}_{a-2}$ commutes with $\gt h_{-1,e}$ and 
acts on $\ff^{a-2}$ via the defining representation. 
%%%% and on $(\ff^{a-2})^*$ as on the dual one. 
The linear subspace $\ff^{a-2}\subset\gt g_{-1,e}$ consists of the vectors
$\xi_i^{1,4}$ with $i>2$, while $(\ff^{a-2})^*\subset\gt g_{-1,e}$ consists 
of the vectors $\xi_1^{i,0}$ with $i>2$. 
Let us choose a basis in $\gt g_{0,e}$ such that the first 
five elements are
$z,\xi_2^{2,0},\xi_1^{1,3},\xi_1^{2,2},\xi_2^{1,4}$ and the other 
$(a-2)^2$ form a basis in $\gt{gl}_{a{-}2}$.
Here  $z,\xi_2^{2,0},\xi_1^{1,3},\xi_1^{2,2},\xi_2^{1,4}$ commute with 
$\ff^{a-2}\oplus(\ff^{a-2})^*$. Note also that 
$\gt{gl}_{a-2}\subset\gt g_{0,e}$ commutes with $\gt h_{-1,e}$.
Hence  
$\rank([\gt g_{0,e},\gt g_{-1,e}])$ is equal 
to the sum of $\rank([\gt h_{0,e},\gt h_{-1,e}])$ 
and the rank of the matrix corresponding to the action 
of $\gt{gl}_{a-2}$ on $\ff^{a-2}{\oplus}(\ff^{a-2})^*$, which is 
$2(a-2)-1$. Summing up, the rank in question is smaller than or equal 
to $1+2(a-2)-1=2(a-2)$. By Lemma~\ref{matrix}, 
$\ind(\gt g_{0,1},\gt g_{-1,e})\ge 4+2(a-2)-2(a-2)=4>3=\ind(\gt g_0,\gt g_1)$.
\end{proof}

Combining Example~\ref{3,0} and Propositions~\ref{3-1}, \ref{m=3;2,a}, 
\ref{m=3-prop} we get the following theorem. 

\begin{thm}\label{m=3}
The only inner automorphisms of $\gt{gl}_n$ of order 
$3$, which have the GIB property correspond to the following 
vectors $\hat r$:
$$
(a,b,0), (a,b,1), (2,2,a).
$$  
\end{thm}

The case $m=3$ is settled now and we pass to higher orders. 
If $m>3$, then there are automorphisms 
of rank $1$ without GIB.

\begin{ex}\label{(2,2,2,1)} 
Suppose that $\hat r=(2,2,2,1)$. Then there is no GIB.
Take a nilpotent element $e\in\gt g_1$ with 
Jordan blocks $(3,3,1)$ and 
$\theta(w_1)=w_1$, $\theta(w_2)=\zeta^2 w_2$, $\theta(w_3)=\zeta w_3$.
Computing the stabiliser we get 
$$
\gt g_{-1,e}=\left< \xi_1^{2,1},\xi_2^{1,1},\xi_2^{3,0}, \xi_3^{2,2} \right>_\ff, 
\qquad
G_{0,e}=(\ff^{^\times})^3\ltimes\exp\left(\left<\xi_1^{2,2},\xi_2^{1,2}\right>_\ff\right).
$$
Note that the nilpotent part of $\gt g_{0,e}$, generated by 
$\xi_1^{2,2}$ and $\xi_2^{1,2}$, acts on 
$\gt g_{-1,e}$ (and hence on its dual) trivially. Therefore 
$\ind(\gt g_{0,e},\gt g_{-1,e})\ge 2>\ind(\gt g_0,\gt g_1)$.
\end{ex}

\begin{ex}\label{(2,2,2,2)}
Suppose that $\hat r=(2,2,2,2)$. Then there is no GIB.
Take a nilpotent element $e\in\gt g_1$ with 
Jordan blocks $(3,3,1,1)$ and 
$\theta(w_1)=w_1$, $\theta(w_2)=\zeta^2 w_2$, $\theta(w_3)=\zeta w_3$,
$\theta(w_4)=\zeta^3 w_4$.
Computing the stabiliser we get 
$$
\gt g_{-1,e}=\left< \xi_1^{2,1},\xi_1^{4,0},\xi_2^{1,1},\xi_2^{3,0}, \xi_3^{2,2},
    \xi_4^{1,2} \right>_\ff,
\qquad 
G_{0,e}=(\ff^{^\times})^4\ltimes\exp\left(\left<\xi_1^{2,2},\xi_2^{1,2}\right>_\ff\right).
$$
The nilpotent part of $\gt g_{0,e}$, generated by $\xi_1^{2,2}$ and $\xi_2^{1,2}$,
acts on  $\gt g_{-1,e}$ (and hence on its dual) trivially. Therefore 
$\ind(\gt g_{0,e},\gt g_{-1,e})\ge 3>\ind(\gt g_0,\gt g_1)$.
\end{ex}

\begin{thm}\label{4-2-1}
Suppose that $m\ge 4$ and either $\rank(G_0,\gt g_1)>1$ or 
$\rank(G_0,\gt g_1)=1$ and $\hat r$ contains a subsequence 
$(a,b,c)$ with $a,b,c\ge 2$. Then the corresponding 
$\theta$-representation has no GIB.
\end{thm}
\begin{proof} 
Suppose first that $\rank(G_0,\gt g_1)>1$. 
Passing to a slice representation (as in Lemma~\ref{slice}), we may assume 
that $\rank(G_0,\gt g_1)=2$. 
Let $e\in\gt g_1$ be a nilpotent 
element with  Jordan blocks $(m-1,m-1,1,\ldots,1)$ 
such that $\theta(w_1)=w_1$, $\theta(w_2)=\zeta^{m-2} w_2$, 
$\theta(w_i)=\zeta w_i$ for $3\le i\le r_1+1$,  and
$\theta(w_j)=\zeta^{m-1} w_j$ for $r_1+2\le j\le r_1+r_{m-1}$.
The remaining generators $w_i$ with $i>r_1+r_2$ 
cannot have eigenvalues $\zeta$ or $\zeta^{m-1}$.
%% with respect to  $\theta$. 
Let $\gt h=\gt{gl}(V)$ be a subalgebra of $\gt g$ 
such that  
$V\subset\VV$ is a $\theta$-invariant subspace of 
dimension $2m+r_1+r_{m-1}-4$,  the restriction of 
$\theta$ to $\gt h$ is an automorphism defined by 
$\hat r_{\gt h}=(2,r_1,2,\ldots,2,r_{m-1})$, and, finally, 
$e\in\gt h_1$. 
Set $a=r_1-1$, $b=r_{m-1}-1$.
We have 
$$
\gt h_{-1,e}=\ff\xi_1^{2,m-3}\oplus\ff\xi_2^{1,1}\oplus(\ff^b{\oplus}(\ff^b)^*)\oplus
   (\ff^a{\oplus}(\ff^a)^*), 
$$
where, for example, $\ff^b$ is generated by $\xi_1^{i,0}$ with $a+3\le i\le r_1+r_{m-1}$.
Suppose that $\xi=\xi_i^{j,s}\in\gt g_{0,e}$ and 
$[\xi,\gt h_{-1,e}]\ne 0$. 
Then $s=0$ and either $i=j\in\{1,2\}$ or 
$\xi\in\gt{gl}_a{\oplus}\gt{gl}_b\subset\gt h_{0,e}$.
In any case $\xi\in\gt h_{0,e}$. 
Hence $\ind(\gt g_{0,e},\gt g_{-1,e})\ge\ind(\gt h_{0,e},\gt h_{-1,e})$.
Let $H\subset G$ be a connected subgroup with $\Lie H=\gt h$.
Then 
$$
H_{0,e}=(\ff^{^\times})^2{\times}\GL_a{\times}\GL_b\ltimes
\exp\left(\left<\xi_1^{2,m-2},\xi_2^{1,2}\right>_\ff\right).
$$
The nilpotent part of $\gt h_{0,e}$ acts on $\gt h_{-1,e}$ trivially. 
Hence 
$\ind(\gt h_{0,e},\gt h_{-1,e})=3$. Therefore 
$\ind(\gt g_{0,e},\gt g_{-1,e})\ge 3>2=\ind(\gt g_0,\gt g_1)$.

Suppose now that $\rank(G_0,\gt g_1)=1$. Without loss of generality we may assume 
that $r_0,r_1,r_2\ge 2$. %%%% and $r_3=1$. 
Take $e\in\gt g_{1}$ with Jordan blocks 
$(m-1,3,1,\ldots,1)$ such that $\theta(w_1)=\zeta^2 w_1$, 
$\theta(w_2)=w_2$, and $\theta(w_i)=\zeta w_i$ for $3\le i\le r_1+1$,
Let $\gt a\subset\gt g_{-1,e}$ be a subspace 
generated by the vectors 
$\xi_1^{2,1}$, $\xi_2^{1,m-3}$, and $\xi_1^{i,0},\xi_i^{1,m-2}$ with 
$3\le i\le r_1+1$. Then $[\gt g_{0,e},\gt a]\subset\gt a$ 
and $G_{0,e}$ acts on $\gt a$ 
as $\ff^{^\times}{\times}\GL_p$ on 
$\ff{\oplus}\ff^*\oplus\ff^a{\oplus}(\ff^a)^*$,
where $a=r_1-1$.
Hence  $\ind(\gt g_{0,e},\gt g_{-1,e})\ge 2>1$. 
\end{proof}

\begin{prop}\label{4-gr}
Combining Theorem~\ref{4-2-1} with Lemma~\ref{fold}, 
we see that  
an  automorphism $\theta$ of order $6$ and rank zero 
defined by $\hat r=(1,a,b,c,1,0)$ with $a,b,c\ge 2$ has no GIB.
\end{prop}

%%% In a way, $1$'s work as separators. 
If $\rank(G_0,\gt g_1)=1$ and $1$ occurs often enough 
among the $r_i$'s, then the $\theta$-representation has GIB.
For example, an automorphism with 
$\hat r=(1,2,1,2,1,2,\ldots,1,2)$ 
always has GIB. 
In order to formalise the statement, we 
need a result in the rank zero case.

\begin{prop}\label{1-groups-0}
Suppose that $\hat r=(1,r_1,\ldots,r_{m-3},1,0)$
has no substrings $r_i,r_{i+1},r_{i+2}$ with all elements 
being larger than $1$. %%%% (here $i+1$ and $i+2$ are considered modulo $m$). 
In other words, if $r_i>1$, then either $r_{i+1}$ or $r_{i+2}$ must be 
$1$ or $0$. Then the corresponding automorphism $\theta$ has GIB.  
\end{prop}
\begin{proof}
We argue by induction on $m$.
In case $m=3$ the statement is obvious, in cases 
$m=4,5$ it was proved in Lemmas~\ref{gl_q}, \ref{rank0-2gr}.
Suppose that $m>5$. Let $i\ge 1$ be the smallest number such that 
$r_i\le 1$. By the assumptions $i\le 3$.
Set $H=\GL_{r_0}{\times}\ldots{\times}\GL_{r_i}$ 
and $\tilde H=\GL_{r_{i+1}}{\times}\ldots{\times}\GL_{r_{m-2}}$.
Then $G_0=H{\times}\tilde H$. As usual $\gt h=\Lie H$ and 
$\tilde{\gt h}=\Lie\tilde H$. 
Let $\gt w\subset\gt g_1$ be the maximal subspace
consisting of $\tilde H$-invariant vectors and 
$\tilde{\gt w}\subset\gt g_1$ its $G_0$-invariant complement.  
We have $\ind(\gt h,\gt w)=0$ and  
the pair $(H,\gt w)$ has GIB by 
the inductive hypothesis. 
%%% Lemma~\ref{gl_q} or 
%%% Lemma~\ref{rank0-2gr} (or by obvious reasons). 
Take $e=x+\tilde x$ with $x\in\gt w$, $\tilde x\in\tilde{\gt w}$. 
According to Lemma~\ref{summa}, there is $y\in\gt w$ 
such that $[\gt h_x,w]+[\gt h,x]=\gt w$.
Let $T,\tilde T$ be the central tori of $H$ and $\tilde H$,
respectively, and set $\gt t=\Lie T$, $\tilde{\gt t}=\Lie\tilde T$.
Note that $T$ acts trivially on $\gt w$.
In particular, $T\subset (H_x)_y$.
The action of $T{\times}\tilde H$ on $\tilde{\gt w}$
is a $\theta$-representation corresponding to 
a vector $(r_i,r_{i+1},\ldots,r_{m-3},1,0)$.
Therefore it has GIB by the inductive hypothesis and 
there is $\tilde y\in\tilde{\gt w}$ such that 
$[(\gt t{\oplus}\tilde{\gt h})_{\tilde x},\tilde y]+[\gt t{\oplus}\tilde{\gt h},\tilde x]=\tilde{\gt w}$.
Combining these two equalities we get 
$$
\begin{array}{l}
[\gt g_0,e]+[\gt g_{0,e},y+\tilde y]=
[\gt g_0,e]+[\gt t{\oplus}\tilde{\gt h},\tilde x]+[\gt g_{0,e},y+\tilde y]+
[(\gt t{\oplus}\tilde{\gt h})_{{\tilde x}},\tilde y]=  \\
\enskip = \tilde{\gt w}+[\gt h,x]+[\gt g_{0,e},y]=  \tilde{\gt w}+[\gt h,x]+[\gt g_{0,e}+\tilde{\gt h},y]. \\
\end{array}
$$
It remains to understand $[\gt g_{0,e}+\gt h,y]$. 
We claim that $\gt h_x\subset\gt g_{0,e}+\gt h$.
In case $\gt h_x$ acts on $\tilde{\gt w}$ trivially, the claim is obvious 
($\gt h_x\subset\gt h_e\subset\gt g_{0,e}$).  If the action is not trivial, then 
$H_x$ acts on $\tilde{\gt w}$ as  $\GL_{r_i}$ and
$[\gt h_x,y]=[\tilde{\gt t},y]$. Therefore 
$\gt h_x\subset(\gt h_x{\oplus}\tilde{\gt t})_y + \tilde{\gt t}\subset
  (\gt g_{0,x})_y+\tilde{\gt t}$. The claim is proved. 

The inclusion  $\gt h_x\subset\gt g_{0,e}+\tilde{\gt h}$ implies that 
$[\gt g_{0,e}+\tilde{\gt h},y]+[\gt h,x]\supset [\gt h_x,y]+[\gt h,x]=\gt w$. 
We see that $[\gt g_0,e]+[\gt g_{0,e},y+\tilde y]=\gt w{\oplus}\tilde{\gt w}=\gt g_1$
and hence $(G_0,\gt g_1)$ has GIB by Lemma~\ref{summa} .
\end{proof}

Combining Proposition~\ref{1-groups-0} with Lemma~\ref{fold},
we get the following.

\begin{prop}\label{1-groups-1}
Suppose that $\hat r=(1,r_1,\ldots,r_{m-2},r_{m-1})$
has no substrings $r_i,r_{i+1},r_{i+2}$ with all elements 
being larger than $1$. 
Then the corresponding automorphism $\theta$ has GIB.
\end{prop}

The rank zero case we leave open. Mainly due to the
following observations. 

\begin{lm}\label{q} Let $q$ and $n$ be natural numbers such that 
$2q\le n$. Then the representation of 
$\gt{gl}_n$ on $V=(q\ff^n)^*$, and hence on $V^*$, has GIB.
\end{lm}
\begin{proof}
Since $q$ is smaller than $n$, the $\GL_n$-orbits 
on $V$ are classified by the matrix rank. 
Let $v\in V$ be a matrix of rank $p$ ($p\le q$). 
In case $p=q$ the orbit is of the maximal dimension 
and it satisfies GIB. Assume $p<q$. 
The quotient space $V/\gt{gl}_n{\cdot}v$ is isomorphic 
to $q\ff^{n-p}$ and $(\GL_n)_v$ acts on it as 
$\GL_{n-p}$. We have $q<(n-p)$, because 
$q+p<2q\le n$. Hence there is an open $(\GL_n)_v$-orbit 
in the quotient $V/\gt{gl}_n{\cdot}v$. 
\end{proof} 

Lemma~\ref{q} leads to the following propagation property. 

\begin{prop}\label{prop}
Suppose that a zero rank automorphism 
$\theta$ of $\gt g=\gt{gl}(\VV)$ with 
$\hat r=(0,r_1,r_2,\ldots,r_{m-1})$ has GIB and 
$q\ge 2r_{m-1}$. Then 
an automorphism of $\gt h=\gt{gl}(\VV{\oplus}\ff^q)$ 
of order $m+1$ defined by 
$\hat r':=(0,r_1,r_2,\ldots,r_{m-1},q)$ 
has GIB as well. 
\end{prop}
\begin{proof}
We have $\gt h_0=\gt g_0\oplus\gt{gl}_q$
and $\gt h_1=\gt g_1\oplus V$, where $V=(\ff^{r_{m-1}})^*{\otimes}\ff^q$.
Take $x\in\gt h_1$. It decomposes as $x=e+v$, 
where $e\in\gt g_1$ and $v\in V$. 
Clearly $(\gt{gl}_q)_v\subset\gt h_{0,x}$.
Because $q\ge r_{m-1}$, we have also
$\gt g_{0,e}\subset\gt h_{0,x}+\gt{gl}_q$.
Both pairs $(G_0,\gt g_1)$ and $(\GL_q,V^*)$
have GIB. Due to Lemma~\ref{summa}, there are
$\xi\in\gt g_1$, $w\in V$ such that 
$[\gt g_{0,e},\xi]+[\gt g_0,e]=\gt g_1$ and 
$(\gt{gl}_q)_v{\cdot}w+\gt{gl}_q{\cdot}v=V$.
Set $y=\xi+w$. 
Then 
$$
\begin{array}{l}
[\gt h_0,x]+[\gt h_{0,x},y]=[\gt{gl}_q,v]+[(\gt{gl}_q)_v,w]+[\gt g_0,x]+[\gt g_{0,x},y]=
V+[\gt g_0,e]+[\gt g_{0,x},\xi]= \\
\enskip =V+[\gt g_0,e]+[(\gt g_{0,x}+\gt{gl}_q),\xi]\supset 
V+[\gt g_0,e]+[\gt g_{0,e},\xi]=\gt h_1. \\
\end{array}
$$ 
Thereby $(H,\gt h_1)$ has GIB by Lemma~\ref{summa}. 
\end{proof}

It seems that in the rank zero case 
the GIB property depends on the entries $r_i$ 
in a rather bizarre way. However computations in small dimensions 
indicate that most likely GIB holds if the order of $\theta$ is $4$ or $5$.  

\begin{cnj} The GIB property holds for rank zero inner automorphisms 
of $\gt{gl}_n$
of orders $4$ and $5$. 
\end{cnj}

\begin{rmk}
Semisimple automorphisms of infinite order yield exactly the same 
$\theta$-representations as finite order automorphisms of rank zero. 
Therefore the problem of checking GIB also remains open for infinite order 
automorphisms. 
\end{rmk}

%%%%%%%%%%%%%%%%%%%%%%%%%%%%%%%%%%%%%%%%%%%%%%%%%%%%%%%%%%%%%%%%%%%%%%%%%%%%
%%%%%%%%%%%%%%%%%%%%%%%%%%%   SECTION %%%%%%%%%%%%%%%%%%%%%%%%%%%%%%%%%%%%%%
\section{GIB in the exceptional types}\label{sec-ex}
%%%%%%%%%%%%%%%%%%%%%%%%%%%%%%%%%%%%%%%%%%%%%%%%%%%%%%%%%%%%%%%%%%%%%%%%%%%%

We have implemented Algorithm~\ref{alg_gib} in {\sf GAP}4, using the
functionality for listing nilpotent orbits of $\theta$-groups present
in the {\sf SLA} package (\cite{sla}). 
For computing the rank of $\widetilde{A}(\underline{a})$, where 
the $a_i$ are indeterminates of a function field, we have used
{\sc Magma} (\cite{magma}). 

We have used this implementation to
find automorphisms $\theta$ of the Lie algebras of exceptional type for which
$\rank(G_0,\g_1) > 0$ and $(G_0,\g_1)$ has GIB. 

In Tables~\ref{tab:gibE6} to \ref{tab:gibG2} the Kac diagrams of the 
automorphisms that we found to have GIB are listed. The explanation of
the content of these tables is as follows.
Since we restrict ourselves to automorphisms of positive rank, 
the labels of the Kac diagrams are
0,1. Hence we give these labels by colouring the nodes in the Kac diagram:
black means that the label is 1, otherwise the label is 0. Note that we
can restrict to automorphisms of order less then the Coxeter number.
Indeed, if the order is equal to that number, then
$G_0$ will be a torus, and hence $(G_0,\g_1)$ has
GIB by \cite[Proposition 1.3]{py-gib}. For higher orders $(G_0,\g_1)$ has
rank zero. In the tables, the first
column has the order of $\theta$, and the second column its Kac diagram.
The third column has the rank of $(G_0,\g_1)$. Moreover, in order to save
space, we put two sets of columns next to each other.

\begin{rmk}
Tables \ref{tab:gibE6} and \ref{tab:gibE7} (inner automorphisms of 
respectively, $\GR{E}{6}$ and $\GR{E}{7}$) contain only one automorphism
of order 2. In \cite{py-gib} it was left open whether or not these cases
have GIB. We conclude that they do.
\end{rmk}

%%\begin{rmk}
In all cases where it was necessary to compute the rank of a matrix
$\widetilde{A}(\underline{a})$, with $a_i$ indeterminates in a function field,
this proved to be a straightforward calculation, except for two cases, both
in $\GR{E}{8}$. In those cases Algorithm \ref{alg_gib} establishes that
GIB does not hold with high probability. Furthermore, exactly one nilpotent
orbit is found that very probably causes GIB to fail. However, it proved to
be a too demanding calculation to compute the rank of  
$\widetilde{A}(\underline{a})$, with $a_i$ indeterminates in a function field.
These two cases are examined in detail in
Examples~\ref{E8-I}, \ref{E8-II}, where we show that they do not have GIB. 
We conclude that the following theorem holds.
%%%%\end{rmk}

\begin{thm}\label{exc-thm}
Tables~\ref{tab:gibE6} to \ref{tab:gibG2} contain the Kac diagrams of all positive
rank automorphisms of the Lie algebras of exceptional type, such that $(G_0,\g_1)$
has GIB. 
\end{thm}

%%%%%% E_8 -- I %%%%%%%%%%%%%%%%%%%%

\begin{ex}\label{E8-I}
Let  $\theta$ be an automorphism of $\GR{E}{8}$ with the following Kac diagram:

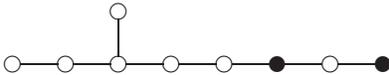
\begin{picture}(142,40)
  \put(3,0){\circle{6}}
  \put(23,0){\circle{6}}
  \put(43,0){\circle{6}}
  \put(63,0){\circle{6}}
  \put(83,0){\circle{6}}
  \put(43,20){\circle{6}}
  \put(103,0){\circle*{6}}
  \put(123,0){\circle{6}}
  \put(143,0){\circle*{6}}

  \put(6,0){\line(1,0){14}}
  \put(26,0){\line(1,0){14}}
  \put(46,0){\line(1,0){14}}
  \put(66,0){\line(1,0){14}}
  \put(43,3){\line(0,1){14}}
\put(86,0){\line(1,0){14}}
\put(106,0){\line(1,0){14}}
\put(126,0){\line(1,0){14}}
\end{picture}

\vskip0.5ex

\noindent
Then $(G_0,\gt g_1)$ does not have GIB. 
By the computations made with Algorithm~\ref{alg_gib},
there is a unique orbit that is possibly bad.  
Here $G_0 = E_6 \times\SL_2 \times T$ with $T\cong\ff^{^\times}$ and
$\gt g_1 = V_1\oplus V_2$ with $V_1=\ff^{27}{\otimes}\ff^2$, $V_2=\ff^2$.
%%%%% Here $\ff^{27}$ is the smallest representation of $E_6$. 
The representation of the group $E_6$ on $\ff^{27}$ is of index $1$
and for generic $v\in\ff^{27}$ the stabiliser $(E_6)_v$ 
is reductive and of type $F_4$, 
%%%% reductive generic stabiliser, isomorphic to $F_4$, 
see e.g. \cite{elashvili}.
%%% the index of this representation is one. 
Let $v\in\ff^{27}$ be such that $(E_6)_v=F_4$.
As is very well known, $F_4$ is the subgroup of stable points for the diagram 
automorphism of $E_6$. Hence its normaliser in $E_6$ 
coincides with $F_4$ up to a connected component. Therefore $\ff v$
is not contained in 
the tangent space $[\GR{E}{6},v]$.
%%%  does not contain ,
%% in other words $\ff v \cap T_v(E_6{\cdot}v) = 0$.   
This implies that $T\times E_6$ acts on $\ff^{27}$ with an 
open orbit. 

Let $t\in T$ and $v_1\in V_1$, $v_2\in V_2$. Then 
$t{\cdot}(v_1+v_2)=tv_1+t^{-3}v_2$.
Two copies of $\ff^2$ are canonically isomorphic as $\SL_2$-modules. 
Using this isomorphism we write 
a representative of the bad orbit as %%%% we take 
$e=v{\otimes}w + w$,
where $w\in\ff^2$ and %% (and it is a kind of the same vector on both places),
$v\in\ff^{27}$ is generic, i.e., such that $(E_6)_v=F_4$.

As we already know, there are no elements $\xi\in\GR{E}{6}$ 
such that $\xi{\cdot}v\ne 0$ and $\xi{\cdot}v\in\ff v$.
Therefore 
$\gt g_{0,v{\otimes}w}=\GR{F}{4}\oplus\gt h\oplus\gt u$, 
where $\gt u\subset\gt{sl}_2$ is the Lie algebra 
of a unipotent subgroup, $\gt h=\ff$, and $\gt h$ 
is embedded diagonally into $\gt t\oplus\gt{sl}_2$  
(here $\gt t=\Lie T$). 
Note that for $h\in\gt h$ and $w\in V_2$ being the same as above, we have
$h{\cdot}w=-4hw$. Hence 
$\gt g_{0,e}=\GR{F}{4}\oplus\gt u$.

The tangent space $[\gt g_0,e]$ is 
equal to 
$$
[\GR{E}{6},v]{\otimes}w\oplus \ff(v{\otimes}w'+w')\oplus \ff v{\otimes}w \oplus \ff w,
$$
where $w'\in\ff^2$ is a vector non-proportional to $w$.
Hence $\gt g_1/[\gt g_0,e]\cong \ff^{26}{\otimes}w' \oplus \ff$, where 
the line $\ff$ is an anti-diagonal in
$\ff(v{\otimes}w') \oplus \ff w'$.
For the nilpotent part $\gt u\subset\gt{sl}_2$ of $\gt g_{0,e}$, 
holds $\gt u{\cdot}w'=w$. Therefore it acts on 
$\gt g_1/[\gt g_0,e]$ trivially.
The representation of $F_4$ on 
the space $\ff^{26}{\otimes}w'\oplus\ff$
is the sum of the trivial representation and the simplest one, which has
index two. Thereby 
$\ind(\gt g_{0,e},\gt g_{-1,e})=3$.  

There are many ways to see that $\ind(\gt g_0,\gt g_1)=2$. 
One of them is to take a slightly modified element in 
$\gt g_1$, namely $x=v{\otimes}w + w'$
with $w$ and $w'$ being linear independent. 
Then $G_{0,x}$-action on $\gt g_1/[\gt g_0,x]$ is the same as the 
action of $F_4$
on $\ff^{26}$ (there is no additional line $\ff$ here). 
We have $\ind(\gt g_{0,x},\gt g_{-1,x})=2$ and since 
the stabiliser $G_{0,x}$ is reductive, 
$\ind(\gt g_{0,x},\gt g_{-1,x})=\ind(\gt g_0,\gt g_1)$, 
see \cite[Proposition~1.1]{py-gib}. 
\end{ex}

%%%%%%%%%%%%%%%%%%%%%%%%%%%%%%%% E_8 -- II %%%%%%%%%%%%%%%%%%%%%%%%%%

\begin{ex}\label{E8-II}
Let $\theta$ be an automorphism of $\GR{E}{8}$ corresponding to the following 
Kac diagram:

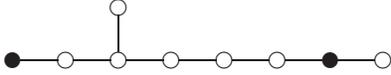
\begin{picture}(142,40)
  \put(3,0){\circle*{6}}
  \put(23,0){\circle{6}}
  \put(43,0){\circle{6}}
  \put(63,0){\circle{6}}
  \put(83,0){\circle{6}}
  \put(43,20){\circle{6}}
  \put(103,0){\circle{6}}
  \put(123,0){\circle*{6}}
  \put(143,0){\circle{6}}

  \put(6,0){\line(1,0){14}}
  \put(26,0){\line(1,0){14}}
  \put(46,0){\line(1,0){14}}
  \put(66,0){\line(1,0){14}}
  \put(43,3){\line(0,1){14}}
\put(86,0){\line(1,0){14}}
\put(106,0){\line(1,0){14}}
\put(126,0){\line(1,0){14}}
\end{picture}

\vskip0.8ex

\noindent
Then $(G_0,\gt g_1)$ does not have GIB.
Here the order of $\theta$ is $4$, $G_0=\Spin_{12}{\times}\SL_2{\times}\ff^{^\times}$,
and $\rank(G_0,\gt g_1)=2$. A
suspicious orbit $G_0{\cdot}e\subset\gt g_1$ was found in accordance with 
Algorithm~\ref{alg_gib}.  It has dimension $29$ and if we include 
$e$ into an $\gt{sl}_2$-triple $\left<e,h,f\right>$ with 
$h\in\gt g_0$, $f\in\gt g_{-1}$, then $h\in\gt{so}_{12}$ and
%%% Consider a nilpotent element $e\in\gt g_1$
%%% with a characteristic 
%%% $h=2h_2+2h_3+4h_4+4h_5+4h_6+2h_7$.  
%%%% in $\gt g_0$. 
%%% The $G_0$-orbit of $e$ is of dimension $29$.
%%%
%%% The characteristic $h$ is an element of $\gt{so}_{12}$ and
the characteristic
$h$ acts in the defining representation $\ff^{12}$ of $\gt{so}_{12}$  
as a semisimple matrix with eigenvalues $(2,2,-2,-2,0^8)$.
Let $\gt t\subset\gt{so}_{12}$ be a maximal torus containing $h$.
We also fix a Borel subalgebra $\gt b\subset\gt{so}_{12}$
containing $\gt t$. 
Replacing $h$ by a $G_0$-conjugate element, if necessary, 
we may (and will) assume that $\esi_1(h)=\esi_2(h)=2$,
$\esi_3(h)=\esi_4(h)=\esi_5(h)=\esi_6(h)=0$ for the standard basis
$\{\esi_1,\ldots,\esi_6\}$ of $\gt t^*$.

We have a $G_0$-invariant decomposition
$\gt g_1=V\oplus W$, where the semisimple part 
of $G_0$ acts on $V=\ff^{32}$ via  
%% As a representation of $G_0$ the subspace 
%%% $\gt g_1$ decomposes into a direct sum of $V=\ff^{32}$, the 
a half-spin representation of $\Spin_{12}$, and on  
$W=\ff^{12}{\otimes} \ff^2$ via the tensor product of 
the defining representations. 
%%% As it makes no difference,
The semisimple element $h$ is invariant under the diagram 
automorphism of $D_6$. Therefore the picture would not 
change if we replace one half-spin representation by another
(this would be just a different choice of the simple roots for $\gt{so}_{12}$). 
It is more convenient to assume that 
%%%% we assume that 
the highest weight $\lambda$ of $V$ 
is equal to $(\esi_1+\esi_2+\ldots+\esi_6)/2$. 
The other weights of $V$ are 
$(\sum\limits_{i=1}^{6}\pm\esi_i)/2$ with even number of minus signs and
each weight space is one-dimensional. 
Let $v_\lambda\in V$ be a highest weight vector. 

In order to identify $e$ in terms of $V$ and $W$,
we need to understand the  
subspace $\gt g_1(2)$, where $2$ stands for 
the eigenvalue of $\ad(h)$. 
Under the action of $G_{0,h}=\Spin_8\times \GL_2 \times \GL_2$
the subspaces $V$ and $W$ decompose as
$V=V_1\oplus V_2\oplus V_3$ and $W=W_1\oplus W_2\oplus W_3$, 
where 
$$
\begin{array}{l}
V_1\cong V_3\cong \ff^8_+ \ \text{   
%%``positive'' half-spin representations of $\Spin_8$, 
%%i.e., with 
with $\ff^8_+$ having the highest weight } \
(\esi_3+\esi_4+\esi_5+\esi_6)/2; \\
V_2\cong\ff^2\otimes\ff^8_- \ \text{ with $\ff^8_-$ having 
the highest weight } \  
(\esi_3+\esi_4+\esi_5-\esi_6)/2; \\ 
W_1\cong W_3\cong \ff^2{\otimes}\ff^2, \enskip W_2\cong\ff^8{\otimes}\ff^2
\ \text{ with $\ff^8$ having the highest weight } \  \esi_3.
\end{array} 
$$
We assume that   
$v_\lambda$ is a highest weight vector in $V_1$ (not in $V_3$).
Note that $(\sum\limits_{i=1}^{6}\pm\esi_i)(h)=4$ if and only 
if $\esi_1$ and $\esi_2$ are taken with the sign $+$.
Therefore $\gt g_1(2)=V_1\oplus W_1$ and
the group $G_{0,h}$
indeed acts on it with an open orbit. (This proves
the existence of a nilpotent element $e\in\gt g_1$ with
the characteristic $h$.) As a representative 
of the open $G_{0,h}$-orbit we choose 
$v+w$, where $v\in V_1$ is the sum 
of $v_\lambda$ and a lowest weight vector  
$v_\mu\in V_1$ (with respect to $\Spin_8$), i.e.,
$\mu=(\esi_1+\esi_2-\esi_3-\esi_4+\esi_5+\esi_6)/2$ 
as a weight of $\Spin_{12}$;
and $w\in W_1(\cong M_{2,2}(\ff^2))$ is the identity matrix.
%%% (Just to check: $(G_{0,h})_e=\SL_2\times \Spin_7$ is 
%%% reductive.) 

Let $L\cong\GL_2$ be a normal subgroup of $G_0$.
By a direct computation we get $G_{0,w}=\GL_2\times \Spin_8 \ltimes N$, where 
$N$ is the unipotent radical of a parabolic 
subgroup  
$P\subset \Spin_{12}$ with the Levi part of type 
$A_1\times D_4$ and $\gt b\subset\Lie P$; 
and the $\GL_2$-factor is embedded diagonally
%% diagonally 
into $P\times L$. 
%% the product of two normal subgroups $\GL_2\lhd P$ and 
%%% Levi of %%% the very same parabolic 
%%% $P$ and the normal factor 
%%% into $\GL_2\lhd G_0$.   
%
%
%
%%%% Clearly $\dim G_0{\cdot}w=\dim N+4=21$ and 
%%%% $\dim (W/[\gt g_0,w])=3$. 
%%%% On the quotient space $W/[\gt g_0,w]$ 
%%%% the stabiliser $G_{0,w}$ acts via the adjoint 
%%%% representation of $\SL_2$. 
%
% We have $\dim(G_0{\cdot}e)=\dim (G_0{\cdot}w)+\dim (G_{0,w}{\cdot}v)$.
% Therefore $\dim (G_{0,w}{\cdot}v)=8$. 
%%%%%  $G_{0,e}=(G_{0,w})_v$ the 
%%% quotient space $\gt g_1/[\gt g_0,e]$ is 
%%% isomorphic to $V/[\gt g_{0,w},v]\oplus W/[\gt g_0,w]$. 
%%
%%
%% Hence one needs to understand the 
%% $G_{0,w}$ orbit of $v$. It is of dimension 
%% $29-21=8$. {\it \ In fact, it coincides with the 
%% orbit of the reductive part. One can show (using weights) that 
We claim that $N{\cdot}v=v$ or, what is the same, that 
$[\gt n,v]=0$ for $\gt n=\Lie N$. 
Because $v_\lambda$ is the highest weight vector, and $\gt n$ is contained 
in the nilpotent radical of $\gt b$, we have $[\gt n,v_\lambda]=0$.
The Lie algebra $\gt n$ consists of weight-spaces 
with weights $\esi_1+\esi_2$ and 
$\esi_1\pm\esi_j$, $\esi_2\pm\esi_j$ with $3\le j\le 6$. 
Clearly each weight of $[\gt n,v_\mu]$ 
has a coefficient $3/2$ in front of $\esi_1$ or $\esi_2$.
Hence $[\gt n,v_\mu]=0$ and $[\gt n,v]=0$. 
%
%  
%
%%% $Nv=v$ or that the stabiliser in the reductive part 
%%% is $SL_2\times Spin_7$.  \  In whatever way we prove it, 
The stabiliser of $v$ in the reductive part of $G_{0,w}$, i.e.,
in $\GL_2\times \Spin_8$, is equal to $\SL_2\times\Spin_7$.
Therefore $G_{0,e}=\SL_2\times \Spin_7\ltimes N$. 

Our final goal is to understand the quotient space 
$\gt g_1/[\gt g_0,e]$. 
Let $N_-$ be the unipotent radical of an opposite to 
$P$ parabolic.  Set $\gt n_-:=\Lie N_-$,
$\gt l:=\Lie L$.  
Then $\gt g_0=\gt g_{0,w}\oplus\gt n_-\oplus\gt l$ and 
$$
[\gt g_0,e]=[\gt g_{0,w},v]\oplus[\gt n_-,e]\oplus[\gt l,w]=
  V_1\oplus[\gt n_-,e]\oplus W_1.
$$
Note that $[\gt n_-,e]=\ff^{16}\oplus\ff$, 
where $\ff^{16}$ is embedded diagonally into $V_2\oplus W_2$,
and $\ff$ is embedded diagonally into $V_3^{\Spin_7}\oplus W_3^{\SL_2}$.
Therefore $\gt g_1/[\gt g_0,e]=(\ff\oplus\ff^3\oplus\ff^7)\oplus\ff^8{\otimes}\ff^2$ 
as an $(\SL_2{\times}\Spin_7)$-module. Here $\SL_2$ acts on 
$\ff^3$ via the adjoint representation.  
Moreover, in the quotient 
$\gt n{\cdot}(\ff{\oplus}\ff^3{\oplus}\ff^7)\subset \ff^8{\otimes}\ff^2$
and $\gt n{\cdot}(\ff^8{\otimes}\ff^2)=0$.   
%%%% 
%
%
% Hence
% $\gt g_1/[\gt g_0,e]\cong V_2\oplus V_3\oplus\gt{sl}_2$. 
% Let $\gt n=\Lie N$ be the Lie algebra of $N$. 
%%% Note that $\gt n=\ff^2{\otimes}\ff^8\oplus\ff$ as 
%%% an $SL_2\times SO_8$-module, where $\ff^8=R(\esi_1)$ 
%%%% is the defining representation of $SO_8$. 
%%%%%% In the quotient $V/V_1$ we have ....
% Then $[\gt n,V_3]\subset V_2\oplus V_1$. 
This implies that 
$\ind(\gt g_{0,e},(\gt g_1/[\gt g_0,e])^*)\ge 
 \ind(\gt{so}_7{\oplus}\gt{sl}_2,\ff{\oplus}\ff^3{\oplus}\ff^7)=3$. 
Since $\ind(\gt g_0,\gt g_1)=2$, we conclude that this
pair $(G_0,\gt g_1)$ does not have GIB.  
\end{ex}  

\begin{longtable}{|c|c|c||c|c|c|}
\caption{Inner automorphisms of $\GR{E}{6}$ such that $(G_0,\g_1)$ has GIB.}\label{tab:gibE6}
\endfirsthead
\hline
\multicolumn{6}{|l|}{\small\slshape GIB automorphisms of $\GR{E}{6}$.} \\
\hline 
\endhead
\hline 
\endfoot
\endlastfoot

\hline

$|\theta|$ & Kac diagram & rk & $|\theta|$ & Kac diagram & rk \\
\hline

2 & 
\begin{picture}(90,55)
  \put(3,5){\circle*{6}}
  \put(23,5){\circle{6}}
  \put(43,5){\circle{6}}
  \put(63,5){\circle{6}}
  \put(83,5){\circle{6}}
  \put(43,25){\circle{6}}
  \put(43,45){\circle*{6}}
  \put(6,5){\line(1,0){14}}
  \put(26,5){\line(1,0){14}}
  \put(46,5){\line(1,0){14}}
  \put(66,5){\line(1,0){14}}
  \put(43,8){\line(0,1){14}}
  \put(43,28){\line(0,1){14}}
\end{picture} & 2 &

3 &
\begin{picture}(90,55)
  \put(3,5){\circle{6}}
  \put(23,5){\circle{6}}
  \put(43,5){\circle{6}}
  \put(63,5){\circle{6}}
  \put(83,5){\circle{6}}
  \put(43,25){\circle*{6}}
  \put(43,45){\circle*{6}}
  \put(6,5){\line(1,0){14}}
  \put(26,5){\line(1,0){14}}
  \put(46,5){\line(1,0){14}}
  \put(66,5){\line(1,0){14}}
  \put(43,8){\line(0,1){14}}
  \put(43,28){\line(0,1){14}}
\end{picture} & 1 \\

3 &
\begin{picture}(90,55)
  \put(3,5){\circle*{6}}
  \put(23,5){\circle{6}}
  \put(43,5){\circle{6}}
  \put(63,5){\circle{6}}
  \put(83,5){\circle*{6}}
  \put(43,25){\circle{6}}
  \put(43,45){\circle*{6}}
  \put(6,5){\line(1,0){14}}
  \put(26,5){\line(1,0){14}}
  \put(46,5){\line(1,0){14}}
  \put(66,5){\line(1,0){14}}
  \put(43,8){\line(0,1){14}}
  \put(43,28){\line(0,1){14}}
\end{picture} & 2 &

4 &
\begin{picture}(90,55)
  \put(3,5){\circle*{6}}
  \put(23,5){\circle{6}}
  \put(43,5){\circle{6}}
  \put(63,5){\circle{6}}
  \put(83,5){\circle{6}}
  \put(43,25){\circle*{6}}
  \put(43,45){\circle*{6}}
  \put(6,5){\line(1,0){14}}
  \put(26,5){\line(1,0){14}}
  \put(46,5){\line(1,0){14}}
  \put(66,5){\line(1,0){14}}
  \put(43,8){\line(0,1){14}}
  \put(43,28){\line(0,1){14}}
\end{picture} & 1 \\

4 &
\begin{picture}(90,55)
  \put(3,5){\circle{6}}
  \put(23,5){\circle*{6}}
  \put(43,5){\circle{6}}
  \put(63,5){\circle{6}}
  \put(83,5){\circle{6}}
  \put(43,25){\circle*{6}}
  \put(43,45){\circle{6}}
  \put(6,5){\line(1,0){14}}
  \put(26,5){\line(1,0){14}}
  \put(46,5){\line(1,0){14}}
  \put(66,5){\line(1,0){14}}
  \put(43,8){\line(0,1){14}}
  \put(43,28){\line(0,1){14}}
\end{picture} & 1 &

5 &
\begin{picture}(90,55)
  \put(3,5){\circle*{6}}
  \put(23,5){\circle{6}}
  \put(43,5){\circle{6}}
  \put(63,5){\circle*{6}}
  \put(83,5){\circle{6}}
  \put(43,25){\circle*{6}}
  \put(43,45){\circle{6}}
  \put(6,5){\line(1,0){14}}
  \put(26,5){\line(1,0){14}}
  \put(46,5){\line(1,0){14}}
  \put(66,5){\line(1,0){14}}
  \put(43,8){\line(0,1){14}}
  \put(43,28){\line(0,1){14}}
\end{picture} & 1 \\

5 &
\begin{picture}(90,55)
  \put(3,5){\circle*{6}}
  \put(23,5){\circle{6}}
  \put(43,5){\circle{6}}
  \put(63,5){\circle*{6}}
  \put(83,5){\circle*{6}}
  \put(43,25){\circle{6}}
  \put(43,45){\circle*{6}}
  \put(6,5){\line(1,0){14}}
  \put(26,5){\line(1,0){14}}
  \put(46,5){\line(1,0){14}}
  \put(66,5){\line(1,0){14}}
  \put(43,8){\line(0,1){14}}
  \put(43,28){\line(0,1){14}}
\end{picture} & 1 &

6 &
\begin{picture}(90,55)
  \put(3,5){\circle*{6}}
  \put(23,5){\circle*{6}}
  \put(43,5){\circle*{6}}
  \put(63,5){\circle{6}}
  \put(83,5){\circle{6}}
  \put(43,25){\circle{6}}
  \put(43,45){\circle{6}}
  \put(6,5){\line(1,0){14}}
  \put(26,5){\line(1,0){14}}
  \put(46,5){\line(1,0){14}}
  \put(66,5){\line(1,0){14}}
  \put(43,8){\line(0,1){14}}
  \put(43,28){\line(0,1){14}}
\end{picture} & 1\\

6 &
\begin{picture}(90,55)
  \put(3,5){\circle{6}}
  \put(23,5){\circle*{6}}
  \put(43,5){\circle{6}}
  \put(63,5){\circle*{6}}
  \put(83,5){\circle{6}}
  \put(43,25){\circle*{6}}
  \put(43,45){\circle{6}}
  \put(6,5){\line(1,0){14}}
  \put(26,5){\line(1,0){14}}
  \put(46,5){\line(1,0){14}}
  \put(66,5){\line(1,0){14}}
  \put(43,8){\line(0,1){14}}
  \put(43,28){\line(0,1){14}}
\end{picture} & 1 &

8 &
\begin{picture}(90,55)
  \put(3,5){\circle*{6}}
  \put(23,5){\circle*{6}}
  \put(43,5){\circle{6}}
  \put(63,5){\circle*{6}}
  \put(83,5){\circle{6}}
  \put(43,25){\circle*{6}}
  \put(43,45){\circle*{6}}
  \put(6,5){\line(1,0){14}}
  \put(26,5){\line(1,0){14}}
  \put(46,5){\line(1,0){14}}
  \put(66,5){\line(1,0){14}}
  \put(43,8){\line(0,1){14}}
  \put(43,28){\line(0,1){14}}
\end{picture} & 1\\

8 &
\begin{picture}(90,55)
  \put(3,5){\circle*{6}}
  \put(23,5){\circle{6}}
  \put(43,5){\circle*{6}}
  \put(63,5){\circle*{6}}
  \put(83,5){\circle*{6}}
  \put(43,25){\circle{6}}
  \put(43,45){\circle*{6}}
  \put(6,5){\line(1,0){14}}
  \put(26,5){\line(1,0){14}}
  \put(46,5){\line(1,0){14}}
  \put(66,5){\line(1,0){14}}
  \put(43,8){\line(0,1){14}}
  \put(43,28){\line(0,1){14}}
\end{picture} & 1 &

& & \\%[1ex]

\hline
\end{longtable}

\begin{longtable}{|c|c|c||c|c|c|}
\caption{Outer automorphisms of $\GR{E}{6}$ such that $(G_0,\g_1)$ has GIB.}\label{tab:gibE6out}
\endfirsthead
\hline
\multicolumn{6}{|l|}{\small\slshape GIB outer automorphisms of $\GR{E}{6}$.} \\
\hline 
\endhead
\hline 
\endfoot
\endlastfoot

\hline

$|\theta|$ & Kac diagram & rk & $|\theta|$ & Kac diagram & rk \\
\hline

2 & 
\begin{picture}(135,15)
  \put(10,5){\circle*{6}}
  \put(40,5){\circle{6}}
  \put(70,5){\circle{6}}
  \put(100,5){\circle{6}}
  \put(130,5){\circle{6}}
  \put(13,5){\line(1,0){24}}
  \put(43,5){\line(1,0){24}}
  \put(72,7){\line(1,0){26}}
  \put(72,3){\line(1,0){26}}
  \put(80,1){\Large $<$}
  \put(103,5){\line(1,0){24}}
\end{picture} & 2 &

4 & 
\begin{picture}(135,15)
  \put(10,5){\circle{6}}
  \put(40,5){\circle*{6}}
  \put(70,5){\circle{6}}
  \put(100,5){\circle{6}}
  \put(130,5){\circle{6}}
  \put(13,5){\line(1,0){24}}
  \put(43,5){\line(1,0){24}}
  \put(72,7){\line(1,0){26}}
  \put(72,3){\line(1,0){26}}
  \put(80,1){\Large $<$}
  \put(103,5){\line(1,0){24}}
\end{picture}
& 1 \\

6 & 
\begin{picture}(135,15)
  \put(10,5){\circle{6}}
  \put(40,5){\circle{6}}
  \put(70,5){\circle{6}}
  \put(100,5){\circle*{6}}
  \put(130,5){\circle*{6}}
  \put(13,5){\line(1,0){24}}
  \put(43,5){\line(1,0){24}}
  \put(72,7){\line(1,0){26}}
  \put(72,3){\line(1,0){26}}
  \put(80,1){\Large $<$}
  \put(103,5){\line(1,0){24}}
\end{picture}
& 1 &

6 & 
\begin{picture}(135,15)
  \put(10,5){\circle{6}}
  \put(40,5){\circle*{6}}
  \put(70,5){\circle{6}}
  \put(100,5){\circle{6}}
  \put(130,5){\circle*{6}}
  \put(13,5){\line(1,0){24}}
  \put(43,5){\line(1,0){24}}
  \put(72,7){\line(1,0){26}}
  \put(72,3){\line(1,0){26}}
  \put(80,1){\Large $<$}
  \put(103,5){\line(1,0){24}}
\end{picture}
& 2\\

6 & 
\begin{picture}(135,15)
  \put(10,5){\circle*{6}}
  \put(40,5){\circle*{6}}
  \put(70,5){\circle{6}}
  \put(100,5){\circle{6}}
  \put(130,5){\circle{6}}
  \put(13,5){\line(1,0){24}}
  \put(43,5){\line(1,0){24}}
  \put(72,7){\line(1,0){26}}
  \put(72,3){\line(1,0){26}}
  \put(80,1){\Large $<$}
  \put(103,5){\line(1,0){24}}
\end{picture}
& 1 &

6 & 
\begin{picture}(135,15)
  \put(10,5){\circle{6}}
  \put(40,5){\circle{6}}
  \put(70,5){\circle*{6}}
  \put(100,5){\circle{6}}
  \put(130,5){\circle{6}}
  \put(13,5){\line(1,0){24}}
  \put(43,5){\line(1,0){24}}
  \put(72,7){\line(1,0){26}}
  \put(72,3){\line(1,0){26}}
  \put(80,1){\Large $<$}
  \put(103,5){\line(1,0){24}}
\end{picture}
& 1 \\

8 & 
\begin{picture}(135,15)
  \put(10,5){\circle*{6}}
  \put(40,5){\circle{6}}
  \put(70,5){\circle{6}}
  \put(100,5){\circle*{6}}
  \put(130,5){\circle*{6}}
  \put(13,5){\line(1,0){24}}
  \put(43,5){\line(1,0){24}}
  \put(72,7){\line(1,0){26}}
  \put(72,3){\line(1,0){26}}
  \put(80,1){\Large $<$}
  \put(103,5){\line(1,0){24}}
\end{picture}
& 1 & 

10 & 
\begin{picture}(135,15)
  \put(10,5){\circle{6}}
  \put(40,5){\circle*{6}}
  \put(70,5){\circle{6}}
  \put(100,5){\circle*{6}}
  \put(130,5){\circle*{6}}
  \put(13,5){\line(1,0){24}}
  \put(43,5){\line(1,0){24}}
  \put(72,7){\line(1,0){26}}
  \put(72,3){\line(1,0){26}}
  \put(80,1){\Large $<$}
  \put(103,5){\line(1,0){24}}
\end{picture}
& 1 \\

10 & 
\begin{picture}(135,15)
  \put(10,5){\circle*{6}}
  \put(40,5){\circle{6}}
  \put(70,5){\circle*{6}}
  \put(100,5){\circle{6}}
  \put(130,5){\circle*{6}}
  \put(13,5){\line(1,0){24}}
  \put(43,5){\line(1,0){24}}
  \put(72,7){\line(1,0){26}}
  \put(72,3){\line(1,0){26}}
  \put(80,1){\Large $<$}
  \put(103,5){\line(1,0){24}}
\end{picture}
& 1 &

12 & 
\begin{picture}(135,15)
  \put(10,5){\circle*{6}}
  \put(40,5){\circle*{6}}
  \put(70,5){\circle{6}}
  \put(100,5){\circle*{6}}
  \put(130,5){\circle*{6}}
  \put(13,5){\line(1,0){24}}
  \put(43,5){\line(1,0){24}}
  \put(72,7){\line(1,0){26}}
  \put(72,3){\line(1,0){26}}
  \put(80,1){\Large $<$}
  \put(103,5){\line(1,0){24}}
\end{picture}
& 1 \\

\hline
\end{longtable}

\begin{longtable}{|c|c|c||c|c|c|}
\caption{Inner automorphisms of $\GR{E}{7}$ such that $(G_0,\g_1)$ has GIB.}\label{tab:gibE7}
\endfirsthead
\hline
\multicolumn{6}{|l|}{\small\slshape GIB automorphisms of $\GR{E}{7}$.} \\
\hline 
\endhead
\hline 
\endfoot
\endlastfoot

\hline

$|\theta|$ & Kac diagram & rk & $|\theta|$ & Kac diagram & rk \\
\hline

2 & 
\begin{picture}(130,45)
  \put(3,5){\circle*{6}}
  \put(23,5){\circle{6}}
  \put(43,5){\circle{6}}
  \put(63,5){\circle{6}}
  \put(83,5){\circle{6}}
  \put(63,25){\circle{6}}
  \put(103,5){\circle{6}}
  \put(123,5){\circle*{6}}

  \put(6,5){\line(1,0){14}}
  \put(26,5){\line(1,0){14}}
  \put(46,5){\line(1,0){14}}
  \put(66,5){\line(1,0){14}}
  \put(63,8){\line(0,1){14}}
\put(86,5){\line(1,0){14}}
\put(106,5){\line(1,0){14}}
\end{picture} & 3 &

3 & 
\begin{picture}(130,40)
  \put(3,5){\circle*{6}}
  \put(23,5){\circle*{6}}
  \put(43,5){\circle{6}}
  \put(63,5){\circle{6}}
  \put(83,5){\circle{6}}
  \put(63,25){\circle{6}}
  \put(103,5){\circle{6}}
  \put(123,5){\circle{6}}

  \put(6,5){\line(1,0){14}}
  \put(26,5){\line(1,0){14}}
  \put(46,5){\line(1,0){14}}
  \put(66,5){\line(1,0){14}}
  \put(63,8){\line(0,1){14}}
\put(86,5){\line(1,0){14}}
\put(106,5){\line(1,0){14}}
\end{picture} & 1 \\

4 & 
\begin{picture}(130,40)
  \put(3,5){\circle{6}}
  \put(23,5){\circle*{6}}
  \put(43,5){\circle{6}}
  \put(63,5){\circle{6}}
  \put(83,5){\circle{6}}
  \put(63,25){\circle*{6}}
  \put(103,5){\circle{6}}
  \put(123,5){\circle{6}}

  \put(6,5){\line(1,0){14}}
  \put(26,5){\line(1,0){14}}
  \put(46,5){\line(1,0){14}}
  \put(66,5){\line(1,0){14}}
  \put(63,8){\line(0,1){14}}
\put(86,5){\line(1,0){14}}
\put(106,5){\line(1,0){14}}
\end{picture} & 1 &

4 & 
\begin{picture}(130,40)
  \put(3,5){\circle*{6}}
  \put(23,5){\circle{6}}
  \put(43,5){\circle{6}}
  \put(63,5){\circle{6}}
  \put(83,5){\circle{6}}
  \put(63,25){\circle{6}}
  \put(103,5){\circle*{6}}
  \put(123,5){\circle*{6}}

  \put(6,5){\line(1,0){14}}
  \put(26,5){\line(1,0){14}}
  \put(46,5){\line(1,0){14}}
  \put(66,5){\line(1,0){14}}
  \put(63,8){\line(0,1){14}}
\put(86,5){\line(1,0){14}}
\put(106,5){\line(1,0){14}}
\end{picture} & 1\\

5 & 
\begin{picture}(130,40)
  \put(3,5){\circle*{6}}
  \put(23,5){\circle{6}}
  \put(43,5){\circle{6}}
  \put(63,5){\circle{6}}
  \put(83,5){\circle{6}}
  \put(63,25){\circle*{6}}
  \put(103,5){\circle*{6}}
  \put(123,5){\circle{6}}

  \put(6,5){\line(1,0){14}}
  \put(26,5){\line(1,0){14}}
  \put(46,5){\line(1,0){14}}
  \put(66,5){\line(1,0){14}}
  \put(63,8){\line(0,1){14}}
\put(86,5){\line(1,0){14}}
\put(106,5){\line(1,0){14}}
\end{picture} & 1 & 

5 & 
\begin{picture}(130,40)
  \put(3,5){\circle*{6}}
  \put(23,5){\circle*{6}}
  \put(43,5){\circle{6}}
  \put(63,5){\circle{6}}
  \put(83,5){\circle{6}}
  \put(63,25){\circle{6}}
  \put(103,5){\circle*{6}}
  \put(123,5){\circle{6}}

  \put(6,5){\line(1,0){14}}
  \put(26,5){\line(1,0){14}}
  \put(46,5){\line(1,0){14}}
  \put(66,5){\line(1,0){14}}
  \put(63,8){\line(0,1){14}}
\put(86,5){\line(1,0){14}}
\put(106,5){\line(1,0){14}}
\end{picture} & 1\\

6 & 
\begin{picture}(130,40)
  \put(3,5){\circle*{6}}
  \put(23,5){\circle*{6}}
  \put(43,5){\circle*{6}}
  \put(63,5){\circle{6}}
  \put(83,5){\circle{6}}
  \put(63,25){\circle{6}}
  \put(103,5){\circle{6}}
  \put(123,5){\circle{6}}

  \put(6,5){\line(1,0){14}}
  \put(26,5){\line(1,0){14}}
  \put(46,5){\line(1,0){14}}
  \put(66,5){\line(1,0){14}}
  \put(63,8){\line(0,1){14}}
\put(86,5){\line(1,0){14}}
\put(106,5){\line(1,0){14}}
\end{picture} & 1 & 

6 & 
\begin{picture}(130,40)
  \put(3,5){\circle{6}}
  \put(23,5){\circle{6}}
  \put(43,5){\circle{6}}
  \put(63,5){\circle*{6}}
  \put(83,5){\circle{6}}
  \put(63,25){\circle*{6}}
  \put(103,5){\circle{6}}
  \put(123,5){\circle{6}}

  \put(6,5){\line(1,0){14}}
  \put(26,5){\line(1,0){14}}
  \put(46,5){\line(1,0){14}}
  \put(66,5){\line(1,0){14}}
  \put(63,8){\line(0,1){14}}
\put(86,5){\line(1,0){14}}
\put(106,5){\line(1,0){14}}
\end{picture} & 1\\

6 & 
\begin{picture}(130,40)
  \put(3,5){\circle{6}}
  \put(23,5){\circle*{6}}
  \put(43,5){\circle{6}}
  \put(63,5){\circle*{6}}
  \put(83,5){\circle{6}}
  \put(63,25){\circle{6}}
  \put(103,5){\circle{6}}
  \put(123,5){\circle{6}}

  \put(6,5){\line(1,0){14}}
  \put(26,5){\line(1,0){14}}
  \put(46,5){\line(1,0){14}}
  \put(66,5){\line(1,0){14}}
  \put(63,8){\line(0,1){14}}
\put(86,5){\line(1,0){14}}
\put(106,5){\line(1,0){14}}
\end{picture} & 1 & 

6 & 
\begin{picture}(130,40)
  \put(3,5){\circle{6}}
  \put(23,5){\circle{6}}
  \put(43,5){\circle*{6}}
  \put(63,5){\circle{6}}
  \put(83,5){\circle*{6}}
  \put(63,25){\circle{6}}
  \put(103,5){\circle{6}}
  \put(123,5){\circle{6}}

  \put(6,5){\line(1,0){14}}
  \put(26,5){\line(1,0){14}}
  \put(46,5){\line(1,0){14}}
  \put(66,5){\line(1,0){14}}
  \put(63,8){\line(0,1){14}}
\put(86,5){\line(1,0){14}}
\put(106,5){\line(1,0){14}}
\end{picture} & 1\\

6 & 
\begin{picture}(130,40)
  \put(3,5){\circle*{6}}
  \put(23,5){\circle*{6}}
  \put(43,5){\circle{6}}
  \put(63,5){\circle{6}}
  \put(83,5){\circle{6}}
  \put(63,25){\circle{6}}
  \put(103,5){\circle*{6}}
  \put(123,5){\circle*{6}}

  \put(6,5){\line(1,0){14}}
  \put(26,5){\line(1,0){14}}
  \put(46,5){\line(1,0){14}}
  \put(66,5){\line(1,0){14}}
  \put(63,8){\line(0,1){14}}
\put(86,5){\line(1,0){14}}
\put(106,5){\line(1,0){14}}
\end{picture} & 1 & 

8 & 
\begin{picture}(130,40)
  \put(3,5){\circle*{6}}
  \put(23,5){\circle{6}}
  \put(43,5){\circle*{6}}
  \put(63,5){\circle{6}}
  \put(83,5){\circle*{6}}
  \put(63,25){\circle{6}}
  \put(103,5){\circle{6}}
  \put(123,5){\circle*{6}}

  \put(6,5){\line(1,0){14}}
  \put(26,5){\line(1,0){14}}
  \put(46,5){\line(1,0){14}}
  \put(66,5){\line(1,0){14}}
  \put(63,8){\line(0,1){14}}
\put(86,5){\line(1,0){14}}
\put(106,5){\line(1,0){14}}
\end{picture} & 1\\

8 & 
\begin{picture}(130,40)
  \put(3,5){\circle*{6}}
  \put(23,5){\circle*{6}}
  \put(43,5){\circle{6}}
  \put(63,5){\circle{6}}
  \put(83,5){\circle*{6}}
  \put(63,25){\circle*{6}}
  \put(103,5){\circle{6}}
  \put(123,5){\circle{6}}

  \put(6,5){\line(1,0){14}}
  \put(26,5){\line(1,0){14}}
  \put(46,5){\line(1,0){14}}
  \put(66,5){\line(1,0){14}}
  \put(63,8){\line(0,1){14}}
\put(86,5){\line(1,0){14}}
\put(106,5){\line(1,0){14}}
\end{picture} & 1 & 

8 & 
\begin{picture}(130,40)
  \put(3,5){\circle*{6}}
  \put(23,5){\circle*{6}}
  \put(43,5){\circle*{6}}
  \put(63,5){\circle{6}}
  \put(83,5){\circle{6}}
  \put(63,25){\circle{6}}
  \put(103,5){\circle*{6}}
  \put(123,5){\circle{6}}

  \put(6,5){\line(1,0){14}}
  \put(26,5){\line(1,0){14}}
  \put(46,5){\line(1,0){14}}
  \put(66,5){\line(1,0){14}}
  \put(63,8){\line(0,1){14}}
\put(86,5){\line(1,0){14}}
\put(106,5){\line(1,0){14}}
\end{picture} & 1\\

8 & 
\begin{picture}(130,40)
  \put(3,5){\circle*{6}}
  \put(23,5){\circle{6}}
  \put(43,5){\circle{6}}
  \put(63,5){\circle*{6}}
  \put(83,5){\circle{6}}
  \put(63,25){\circle*{6}}
  \put(103,5){\circle{6}}
  \put(123,5){\circle*{6}}

  \put(6,5){\line(1,0){14}}
  \put(26,5){\line(1,0){14}}
  \put(46,5){\line(1,0){14}}
  \put(66,5){\line(1,0){14}}
  \put(63,8){\line(0,1){14}}
\put(86,5){\line(1,0){14}}
\put(106,5){\line(1,0){14}}
\end{picture} & 1 & 

8 & 
\begin{picture}(130,40)
  \put(3,5){\circle*{6}}
  \put(23,5){\circle*{6}}
  \put(43,5){\circle{6}}
  \put(63,5){\circle{6}}
  \put(83,5){\circle{6}}
  \put(63,25){\circle*{6}}
  \put(103,5){\circle*{6}}
  \put(123,5){\circle*{6}}

  \put(6,5){\line(1,0){14}}
  \put(26,5){\line(1,0){14}}
  \put(46,5){\line(1,0){14}}
  \put(66,5){\line(1,0){14}}
  \put(63,8){\line(0,1){14}}
\put(86,5){\line(1,0){14}}
\put(106,5){\line(1,0){14}}
\end{picture} & 1\\

10 & 
\begin{picture}(130,40)
  \put(3,5){\circle{6}}
  \put(23,5){\circle*{6}}
  \put(43,5){\circle{6}}
  \put(63,5){\circle*{6}}
  \put(83,5){\circle{6}}
  \put(63,25){\circle*{6}}
  \put(103,5){\circle*{6}}
  \put(123,5){\circle{6}}

  \put(6,5){\line(1,0){14}}
  \put(26,5){\line(1,0){14}}
  \put(46,5){\line(1,0){14}}
  \put(66,5){\line(1,0){14}}
  \put(63,8){\line(0,1){14}}
\put(86,5){\line(1,0){14}}
\put(106,5){\line(1,0){14}}
\end{picture} & 1 & 

12 & 
\begin{picture}(130,40)
  \put(3,5){\circle*{6}}
  \put(23,5){\circle*{6}}
  \put(43,5){\circle*{6}}
  \put(63,5){\circle*{6}}
  \put(83,5){\circle{6}}
  \put(63,25){\circle{6}}
  \put(103,5){\circle*{6}}
  \put(123,5){\circle{6}}

  \put(6,5){\line(1,0){14}}
  \put(26,5){\line(1,0){14}}
  \put(46,5){\line(1,0){14}}
  \put(66,5){\line(1,0){14}}
  \put(63,8){\line(0,1){14}}
\put(86,5){\line(1,0){14}}
\put(106,5){\line(1,0){14}}
\end{picture} & 1\\

12 & 
\begin{picture}(130,40)
  \put(3,5){\circle*{6}}
  \put(23,5){\circle*{6}}
  \put(43,5){\circle{6}}
  \put(63,5){\circle*{6}}
  \put(83,5){\circle{6}}
  \put(63,25){\circle*{6}}
  \put(103,5){\circle*{6}}
  \put(123,5){\circle*{6}}

  \put(6,5){\line(1,0){14}}
  \put(26,5){\line(1,0){14}}
  \put(46,5){\line(1,0){14}}
  \put(66,5){\line(1,0){14}}
  \put(63,8){\line(0,1){14}}
\put(86,5){\line(1,0){14}}
\put(106,5){\line(1,0){14}}
\end{picture} & 1 &
&& \\

\hline
\end{longtable}

\begin{longtable}{|c|c|c||c|c|c|}
\caption{Automorphisms of $\GR{E}{8}$ such that $(G_0,\g_1)$ has GIB.}\label{tab:gibE8}
\endfirsthead
\hline
\multicolumn{6}{|l|}{\small\slshape GIB automorphisms of $\GR{E}{8}$.} \\
\hline 
\endhead
\hline 
\endfoot
\endlastfoot

\hline

$|\theta|$ & Kac diagram & rk & $|\theta|$ & Kac diagram & rk \\
\hline

3 & 
\begin{picture}(147,40)
  \put(3,5){\circle{6}}
  \put(23,5){\circle{6}}
  \put(43,5){\circle{6}}
  \put(63,5){\circle{6}}
  \put(83,5){\circle{6}}
  \put(43,25){\circle{6}}
  \put(103,5){\circle{6}}
  \put(123,5){\circle*{6}}
  \put(143,5){\circle*{6}}

  \put(6,5){\line(1,0){14}}
  \put(26,5){\line(1,0){14}}
  \put(46,5){\line(1,0){14}}
  \put(66,5){\line(1,0){14}}
  \put(43,8){\line(0,1){14}}
\put(86,5){\line(1,0){14}}
\put(106,5){\line(1,0){14}}
\put(126,5){\line(1,0){14}}
\end{picture} & 1 &

5 & 
\begin{picture}(147,40)
  \put(3,5){\circle*{6}}
  \put(23,5){\circle{6}}
  \put(43,5){\circle{6}}
  \put(63,5){\circle{6}}
  \put(83,5){\circle{6}}
  \put(43,25){\circle{6}}
  \put(103,5){\circle{6}}
  \put(123,5){\circle*{6}}
  \put(143,5){\circle*{6}}

  \put(6,5){\line(1,0){14}}
  \put(26,5){\line(1,0){14}}
  \put(46,5){\line(1,0){14}}
  \put(66,5){\line(1,0){14}}
  \put(43,8){\line(0,1){14}}
\put(86,5){\line(1,0){14}}
\put(106,5){\line(1,0){14}}
\put(126,5){\line(1,0){14}}
\end{picture} & 1 \\

6 & 
\begin{picture}(147,40)
  \put(3,5){\circle*{6}}
  \put(23,5){\circle*{6}}
  \put(43,5){\circle{6}}
  \put(63,5){\circle{6}}
  \put(83,5){\circle{6}}
  \put(43,25){\circle{6}}
  \put(103,5){\circle{6}}
  \put(123,5){\circle{6}}
  \put(143,5){\circle{6}}

  \put(6,5){\line(1,0){14}}
  \put(26,5){\line(1,0){14}}
  \put(46,5){\line(1,0){14}}
  \put(66,5){\line(1,0){14}}
  \put(43,8){\line(0,1){14}}
\put(86,5){\line(1,0){14}}
\put(106,5){\line(1,0){14}}
\put(126,5){\line(1,0){14}}
\end{picture} & 1 &

6 & 
\begin{picture}(147,40)
  \put(3,5){\circle{6}}
  \put(23,5){\circle{6}}
  \put(43,5){\circle*{6}}
  \put(63,5){\circle{6}}
  \put(83,5){\circle{6}}
  \put(43,25){\circle{6}}
  \put(103,5){\circle{6}}
  \put(123,5){\circle{6}}
  \put(143,5){\circle{6}}

  \put(6,5){\line(1,0){14}}
  \put(26,5){\line(1,0){14}}
  \put(46,5){\line(1,0){14}}
  \put(66,5){\line(1,0){14}}
  \put(43,8){\line(0,1){14}}
\put(86,5){\line(1,0){14}}
\put(106,5){\line(1,0){14}}
\put(126,5){\line(1,0){14}}
\end{picture} & 1 \\

6 & 
\begin{picture}(147,40)
  \put(3,5){\circle{6}}
  \put(23,5){\circle{6}}
  \put(43,5){\circle{6}}
  \put(63,5){\circle{6}}
  \put(83,5){\circle*{6}}
  \put(43,25){\circle{6}}
  \put(103,5){\circle{6}}
  \put(123,5){\circle*{6}}
  \put(143,5){\circle{6}}

  \put(6,5){\line(1,0){14}}
  \put(26,5){\line(1,0){14}}
  \put(46,5){\line(1,0){14}}
  \put(66,5){\line(1,0){14}}
  \put(43,8){\line(0,1){14}}
\put(86,5){\line(1,0){14}}
\put(106,5){\line(1,0){14}}
\put(126,5){\line(1,0){14}}
\end{picture} & 1 &

6 & 
\begin{picture}(147,40)
  \put(3,5){\circle{6}}
  \put(23,5){\circle{6}}
  \put(43,5){\circle{6}}
  \put(63,5){\circle{6}}
  \put(83,5){\circle{6}}
  \put(43,25){\circle{6}}
  \put(103,5){\circle*{6}}
  \put(123,5){\circle*{6}}
  \put(143,5){\circle*{6}}

  \put(6,5){\line(1,0){14}}
  \put(26,5){\line(1,0){14}}
  \put(46,5){\line(1,0){14}}
  \put(66,5){\line(1,0){14}}
  \put(43,8){\line(0,1){14}}
\put(86,5){\line(1,0){14}}
\put(106,5){\line(1,0){14}}
\put(126,5){\line(1,0){14}}
\end{picture} & 1 \\

8 & 
\begin{picture}(147,40)
  \put(3,5){\circle{6}}
  \put(23,5){\circle{6}}
  \put(43,5){\circle{6}}
  \put(63,5){\circle{6}}
  \put(83,5){\circle*{6}}
  \put(43,25){\circle*{6}}
  \put(103,5){\circle{6}}
  \put(123,5){\circle{6}}
  \put(143,5){\circle*{6}}

  \put(6,5){\line(1,0){14}}
  \put(26,5){\line(1,0){14}}
  \put(46,5){\line(1,0){14}}
  \put(66,5){\line(1,0){14}}
  \put(43,8){\line(0,1){14}}
\put(86,5){\line(1,0){14}}
\put(106,5){\line(1,0){14}}
\put(126,5){\line(1,0){14}}
\end{picture} & 1 &

8 & 
\begin{picture}(147,40)
  \put(3,5){\circle*{6}}
  \put(23,5){\circle*{6}}
  \put(43,5){\circle{6}}
  \put(63,5){\circle{6}}
  \put(83,5){\circle{6}}
  \put(43,25){\circle{6}}
  \put(103,5){\circle{6}}
  \put(123,5){\circle*{6}}
  \put(143,5){\circle{6}}

  \put(6,5){\line(1,0){14}}
  \put(26,5){\line(1,0){14}}
  \put(46,5){\line(1,0){14}}
  \put(66,5){\line(1,0){14}}
  \put(43,8){\line(0,1){14}}
\put(86,5){\line(1,0){14}}
\put(106,5){\line(1,0){14}}
\put(126,5){\line(1,0){14}}
\end{picture} & 1 \\

8 & 
\begin{picture}(147,40)
  \put(3,5){\circle*{6}}
  \put(23,5){\circle{6}}
  \put(43,5){\circle{6}}
  \put(63,5){\circle{6}}
  \put(83,5){\circle{6}}
  \put(43,25){\circle{6}}
  \put(103,5){\circle*{6}}
  \put(123,5){\circle*{6}}
  \put(143,5){\circle*{6}}

  \put(6,5){\line(1,0){14}}
  \put(26,5){\line(1,0){14}}
  \put(46,5){\line(1,0){14}}
  \put(66,5){\line(1,0){14}}
  \put(43,8){\line(0,1){14}}
\put(86,5){\line(1,0){14}}
\put(106,5){\line(1,0){14}}
\put(126,5){\line(1,0){14}}
\end{picture} & 1 &

8 & 
\begin{picture}(147,40)
  \put(3,5){\circle*{6}}
  \put(23,5){\circle{6}}
  \put(43,5){\circle{6}}
  \put(63,5){\circle{6}}
  \put(83,5){\circle{6}}
  \put(43,25){\circle*{6}}
  \put(103,5){\circle{6}}
  \put(123,5){\circle*{6}}
  \put(143,5){\circle*{6}}

  \put(6,5){\line(1,0){14}}
  \put(26,5){\line(1,0){14}}
  \put(46,5){\line(1,0){14}}
  \put(66,5){\line(1,0){14}}
  \put(43,8){\line(0,1){14}}
\put(86,5){\line(1,0){14}}
\put(106,5){\line(1,0){14}}
\put(126,5){\line(1,0){14}}
\end{picture} & 1 \\

10 & 
\begin{picture}(147,40)
  \put(3,5){\circle*{6}}
  \put(23,5){\circle*{6}}
  \put(43,5){\circle{6}}
  \put(63,5){\circle{6}}
  \put(83,5){\circle*{6}}
  \put(43,25){\circle{6}}
  \put(103,5){\circle{6}}
  \put(123,5){\circle{6}}
  \put(143,5){\circle{6}}

  \put(6,5){\line(1,0){14}}
  \put(26,5){\line(1,0){14}}
  \put(46,5){\line(1,0){14}}
  \put(66,5){\line(1,0){14}}
  \put(43,8){\line(0,1){14}}
\put(86,5){\line(1,0){14}}
\put(106,5){\line(1,0){14}}
\put(126,5){\line(1,0){14}}
\end{picture} & 1 &

12 & 
\begin{picture}(147,40)
  \put(3,5){\circle{6}}
  \put(23,5){\circle{6}}
  \put(43,5){\circle*{6}}
  \put(63,5){\circle{6}}
  \put(83,5){\circle{6}}
  \put(43,25){\circle*{6}}
  \put(103,5){\circle*{6}}
  \put(123,5){\circle{6}}
  \put(143,5){\circle{6}}

  \put(6,5){\line(1,0){14}}
  \put(26,5){\line(1,0){14}}
  \put(46,5){\line(1,0){14}}
  \put(66,5){\line(1,0){14}}
  \put(43,8){\line(0,1){14}}
\put(86,5){\line(1,0){14}}
\put(106,5){\line(1,0){14}}
\put(126,5){\line(1,0){14}}
\end{picture} & 1 \\

12 & 
\begin{picture}(147,40)
  \put(3,5){\circle*{6}}
  \put(23,5){\circle*{6}}
  \put(43,5){\circle{6}}
  \put(63,5){\circle{6}}
  \put(83,5){\circle*{6}}
  \put(43,25){\circle{6}}
  \put(103,5){\circle{6}}
  \put(123,5){\circle*{6}}
  \put(143,5){\circle{6}}

  \put(6,5){\line(1,0){14}}
  \put(26,5){\line(1,0){14}}
  \put(46,5){\line(1,0){14}}
  \put(66,5){\line(1,0){14}}
  \put(43,8){\line(0,1){14}}
\put(86,5){\line(1,0){14}}
\put(106,5){\line(1,0){14}}
\put(126,5){\line(1,0){14}}
\end{picture} & 1 &

12 & 
\begin{picture}(147,40)
  \put(3,5){\circle{6}}
  \put(23,5){\circle*{6}}
  \put(43,5){\circle{6}}
  \put(63,5){\circle*{6}}
  \put(83,5){\circle{6}}
  \put(43,25){\circle{6}}
  \put(103,5){\circle{6}}
  \put(123,5){\circle*{6}}
  \put(143,5){\circle*{6}}

  \put(6,5){\line(1,0){14}}
  \put(26,5){\line(1,0){14}}
  \put(46,5){\line(1,0){14}}
  \put(66,5){\line(1,0){14}}
  \put(43,8){\line(0,1){14}}
\put(86,5){\line(1,0){14}}
\put(106,5){\line(1,0){14}}
\put(126,5){\line(1,0){14}}
\end{picture} & 1 \\

12 & 
\begin{picture}(147,40)
  \put(3,5){\circle{6}}
  \put(23,5){\circle{6}}
  \put(43,5){\circle*{6}}
  \put(63,5){\circle{6}}
  \put(83,5){\circle*{6}}
  \put(43,25){\circle{6}}
  \put(103,5){\circle{6}}
  \put(123,5){\circle*{6}}
  \put(143,5){\circle{6}}

  \put(6,5){\line(1,0){14}}
  \put(26,5){\line(1,0){14}}
  \put(46,5){\line(1,0){14}}
  \put(66,5){\line(1,0){14}}
  \put(43,8){\line(0,1){14}}
\put(86,5){\line(1,0){14}}
\put(106,5){\line(1,0){14}}
\put(126,5){\line(1,0){14}}
\end{picture} & 1 &

12 & 
\begin{picture}(147,40)
  \put(3,5){\circle*{6}}
  \put(23,5){\circle{6}}
  \put(43,5){\circle{6}}
  \put(63,5){\circle{6}}
  \put(83,5){\circle*{6}}
  \put(43,25){\circle{6}}
  \put(103,5){\circle*{6}}
  \put(123,5){\circle*{6}}
  \put(143,5){\circle*{6}}

  \put(6,5){\line(1,0){14}}
  \put(26,5){\line(1,0){14}}
  \put(46,5){\line(1,0){14}}
  \put(66,5){\line(1,0){14}}
  \put(43,8){\line(0,1){14}}
\put(86,5){\line(1,0){14}}
\put(106,5){\line(1,0){14}}
\put(126,5){\line(1,0){14}}
\end{picture} & 1 \\

18 & 
\begin{picture}(147,40)
  \put(3,5){\circle*{6}}
  \put(23,5){\circle*{6}}
  \put(43,5){\circle*{6}}
  \put(63,5){\circle{6}}
  \put(83,5){\circle*{6}}
  \put(43,25){\circle{6}}
  \put(103,5){\circle{6}}
  \put(123,5){\circle*{6}}
  \put(143,5){\circle{6}}

  \put(6,5){\line(1,0){14}}
  \put(26,5){\line(1,0){14}}
  \put(46,5){\line(1,0){14}}
  \put(66,5){\line(1,0){14}}
  \put(43,8){\line(0,1){14}}
\put(86,5){\line(1,0){14}}
\put(106,5){\line(1,0){14}}
\put(126,5){\line(1,0){14}}
\end{picture} & 1 & & & \\

\hline
\end{longtable}

\begin{longtable}{|c|c|c||c|c|c|}
\caption{Automorphisms of $\GR{F}{4}$ such that $(G_0,\g_1)$ has GIB.}\label{tab:gibF4}
\endfirsthead
\hline
\multicolumn{6}{|l|}{\small\slshape GIB automorphisms of $\GR{F}{4}$.} \\
\hline 
\endhead
\hline 
\endfoot
\endlastfoot

\hline

$|\theta|$ & Kac diagram & rk & $|\theta|$ & Kac diagram & rk \\
\hline

2 & 
\begin{picture}(135,15)
  \put(10,5){\circle{6}}
  \put(40,5){\circle{6}}
  \put(70,5){\circle{6}}
  \put(100,5){\circle{6}}
  \put(130,5){\circle*{6}}
  \put(13,5){\line(1,0){24}}
  \put(43,5){\line(1,0){24}}
  \put(72,7){\line(1,0){26}}
  \put(72,3){\line(1,0){26}}
  \put(80,1){\Large $>$}
  \put(103,5){\line(1,0){24}}
\end{picture}
& 1 &

3 & 
\begin{picture}(135,15)
  \put(10,5){\circle*{6}}
  \put(40,5){\circle{6}}
  \put(70,5){\circle{6}}
  \put(100,5){\circle{6}}
  \put(130,5){\circle*{6}}
  \put(13,5){\line(1,0){24}}
  \put(43,5){\line(1,0){24}}
  \put(72,7){\line(1,0){26}}
  \put(72,3){\line(1,0){26}}
  \put(80,1){\Large $>$}
  \put(103,5){\line(1,0){24}}
\end{picture}
& 1 \\

3 & 
\begin{picture}(135,15)
  \put(10,5){\circle*{6}}
  \put(40,5){\circle*{6}}
  \put(70,5){\circle{6}}
  \put(100,5){\circle{6}}
  \put(130,5){\circle{6}}
  \put(13,5){\line(1,0){24}}
  \put(43,5){\line(1,0){24}}
  \put(72,7){\line(1,0){26}}
  \put(72,3){\line(1,0){26}}
  \put(80,1){\Large $>$}
  \put(103,5){\line(1,0){24}}
\end{picture}
& 1 &

6 & 
\begin{picture}(135,15)
  \put(10,5){\circle{6}}
  \put(40,5){\circle*{6}}
  \put(70,5){\circle{6}}
  \put(100,5){\circle*{6}}
  \put(130,5){\circle{6}}
  \put(13,5){\line(1,0){24}}
  \put(43,5){\line(1,0){24}}
  \put(72,7){\line(1,0){26}}
  \put(72,3){\line(1,0){26}}
  \put(80,1){\Large $>$}
  \put(103,5){\line(1,0){24}}
\end{picture}
& 1 \\

6 & 
\begin{picture}(135,15)
  \put(10,5){\circle*{6}}
  \put(40,5){\circle*{6}}
  \put(70,5){\circle*{6}}
  \put(100,5){\circle{6}}
  \put(130,5){\circle{6}}
  \put(13,5){\line(1,0){24}}
  \put(43,5){\line(1,0){24}}
  \put(72,7){\line(1,0){26}}
  \put(72,3){\line(1,0){26}}
  \put(80,1){\Large $>$}
  \put(103,5){\line(1,0){24}}
\end{picture}
& 1 &

8 & 
\begin{picture}(135,15)
  \put(10,5){\circle*{6}}
  \put(40,5){\circle*{6}}
  \put(70,5){\circle*{6}}
  \put(100,5){\circle{6}}
  \put(130,5){\circle*{6}}
  \put(13,5){\line(1,0){24}}
  \put(43,5){\line(1,0){24}}
  \put(72,7){\line(1,0){26}}
  \put(72,3){\line(1,0){26}}
  \put(80,1){\Large $>$}
  \put(103,5){\line(1,0){24}}
\end{picture}
& 1 \\

\hline
\end{longtable}

\begin{longtable}{|c|c|c|}
\caption{Automorphisms of $\GR{G}{2}$ such that $(G_0,\g_1)$ has GIB.}\label{tab:gibG2}
\endfirsthead
\hline
\multicolumn{3}{|l|}{\small\slshape GIB automorphisms of $\GR{G}{2}$.} \\
\hline 
\endhead
\hline 
\endfoot
\endlastfoot

\hline

$|\theta|$ & Kac diagram & rk \\
\hline

3 &
\begin{picture}(70,17)
  \put(5,5){\circle*{6}}
  \put(35,5){\circle*{6}}
  \put(65,5){\circle{6}}
  \put(8,5){\line(1,0){24}}
  \put(35,2){\line(1,0){30}}
  \put(38,5){\line(1,0){24}}
  \put(35,8){\line(1,0){30}}
  \put(45,1){\Large $>$}
\end{picture}
& 1 \\

\hline
\end{longtable}


\begin{thebibliography}{33}

\bibitem{magma}
W.~Bosma, J.~Cannon, and C.~Playoust.
\newblock The {M}agma algebra system. {I}. {T}he user language.
\newblock {\em J. Symbolic Comput.}, 24(3-4):235--265, 1997.
\newblock Computational algebra and number theory (London, 1993).

\bibitem{ch-m}
J.-Y.~Charbonnel and A.~Moreau.
The index of centralizers of elements of reductive Lie algebras. {\tt arXiv:0904.1778v1}[math.RT], 2009. 

\bibitem{elashvili} 
A.G.~{\`E}lashvili.
Canonical form and stationary subalgebras of points in general 
position for simple linear Lie groups.
{\it Funkts. Anal. i Prilozh.} 
{\bf 6}\,(1972), no.~1, 51--62; english translation in
{\it Funct. Anal. Appl.} {\bf 6}\,(1972), 44--53. 

\bibitem{gap4}
The GAP~Group.
\newblock {\em {GAP -- Groups, Algorithms, and Programming, Version 4.4}},
  2004. \\
\newblock \verb+(http://www.gap-system.org)+.

%\bibitem{gra6}
%W.~A.~de Graaf.
%\newblock {\em Lie Algebras: Theory and Algorithms}, volume~56 of {\em
%  North-Holland Mathematical Library}.
%\newblock Elsevier Science, 2000.

\bibitem{graaf}
Willem A.~de Graaf.
Computing with nilpotent orbits in simple Lie algebras of exceptional type.
{\it  LMS J. Comput. Math.}, {\bf 11}\,(2008), 280--297. 

\bibitem{sla}
W.~A.~de Graaf.
\newblock {\sf SLA} - computing with {S}imple {L}ie {A}lgebras.
\newblock a {\sf GAP} package, 2009. \\
\newblock \verb+(http://science.unitn.it/~degraaf/sla.html)+.

\bibitem{gra15}
Willem A.~de Graaf.
\newblock Computing representatives of nilpotent orbits of $\theta$-groups.
\newblock %%preprint, 
\\  {\tt arXiv:0905.3149v1}[math.RT], 2009.

\bibitem{helgason}
Sigurdur Helgason.
\newblock {\em Differential geometry, {L}ie groups, and symmetric spaces},
  volume~80 of {\em Pure and Applied Mathematics}.
\newblock Academic Press Inc. [Harcourt Brace Jovanovich Publishers], New York,
  1978.

%\bibitem{hum}
%J.~E. Humphreys.
%\newblock {\em {Introduction to Lie Algebras and Representation Theory}}.
%\newblock Springer Verlag, New York, Heidelberg, Berlin, 1972.

\bibitem{kac_autom}
V.~G. Kac.
\newblock Automorphisms of finite order of semisimple {L}ie algebras.
\newblock {\em Funkcional. Anal. i Prilo\v zen.}, 3(3):94--96, 1969.

\bibitem{Gisela-dis}
G.~Kempken.
Eine Darstellung des K{\"o}chers $\tilde A_{K}$.
Dissertation,  
Rheinischen Friedrich-Wilhelms-Universit\"at Bonn, 1981;  
Bonner Mathematische Schriften {\bf 137}\,(1982).

\bibitem{kr}
B.~Kostant and S.~Rallis.
\newblock Orbits and representations associated with symmetric spaces.
\newblock {\em Amer. J. Math.}, 93:753--809, 1971.

\bibitem{Dima2} 
D.~Panyushev.
The index of a Lie algebra, the centralizer of a nilpotent 
element, and the normalizer of the centralizer. 
{\it Math. Proc. Cambr. Phil. Soc.}, 
(2003)\,{\bf 134}, no.1, 41--59. 

\bibitem{py-gib}
Dmitri~I. Panyushev and Oksana~S. Yakimova.
\newblock The index of representations associated with stabilisers.
\newblock {\em J. Algebra}, 302(1):280--304, 2006.

\bibitem{vin-nilp}
{\`E}.~B. Vinberg.
The classification of nilpotent elements of graded Lie algebras. 
{\it Dokl. Akad. Nauk SSSR}, 225(4):745–748, 1975. 

\bibitem{vinberg}
{\`E}.~B. Vinberg.
\newblock The {W}eyl group of a graded {L}ie algebra.
\newblock {\em Izv. Akad. Nauk SSSR Ser. Mat.}, 40(3):488--526, 709, 1976.
\newblock English translation: Math. USSR-Izv. 10, 463--495 (1976).

 
% \bibitem{cm}
% {\sc D.\,Collingwood} and {W.\,McGovern}, ``Nilpotent orbits in
% semisimple Lie algebras'', Mathematics Series, Van Nostrand Reinhold, 1993.


% \bibitem{ja}
% {\sc J.C.\,Jantzen}, Nilpotent orbits in representation theory, in:
% B.\,Orsted (ed.), ``Representation and Lie theory'', Progr. in
% Math., {\bf 228}, 1--211, Birkh\"auser, Boston 2004.

%%% \bibitem{kr}
% {\sc B.\,Kostant} and {\sc S.\,Rallis}, Orbits and
% representations associated with symmetric spaces,
% {\it Amer. J. Math.} {\bf 93}\,(1971), 753--809.

% \bibitem{py-gib} 
% {\sc D.I.~Panyushev} and  {O.S.~Yakimova},
% The index of representations associated
% with stabilisers, {\it J. Algebra},
% {\bf 302}\,(2006), 280–-304.

%%% \bibitem{vinberg}
%% {\sc {\'E}.B.\,Vinberg}, The Weyl group of a graded Lie algebra, 
%% {\it Izv. Akad. Nauk SSSR Ser.
%% Mat.,} {\bf 40}(3):488--526, 709, 1976. 
%%% English translation: Math. USSR-Izv. 10, 463--495 (1976).

\bibitem{fan}
{\rus O.~Yakimova}. {\rus Indeks centralizatorov {e1}lementov v
klassicheskikh algebrakh Li}. {\rusi Funkc. analiz i ego prilozh.},
{\bf 40}, {\rus N0}\,1 (2006), 52--64 (Russian). English
translation: O.~Yakimova. The centralisers of nilpotent
elements in classical Lie algebras. {\it Funct. Anal. Appl.}, 
{\bf 40} (2006), 42--51.

\end{thebibliography}
\end{document}